\documentclass[11pt]{amsart}

\usepackage{amsmath,amssymb,amsthm}
\usepackage{graphicx}
\usepackage{booktabs}
\usepackage{hyperref}
\usepackage[margin=1in]{geometry}
\usepackage{caption}
\usepackage{enumitem}
\usepackage{xcolor}
\hypersetup{colorlinks=true,linkcolor=blue!70!black,citecolor=green!50!black,urlcolor=blue!70!black}

\theoremstyle{remark}
\newtheorem{remark}{Remark}[section]

\makeatletter
\@ifundefined{Sha}{%
  \DeclareFontFamily{U}{wncy}{}
  \DeclareFontShape{U}{wncy}{m}{n}{<->wncyr10}{}
  \DeclareSymbolFont{cyrletters}{U}{wncy}{m}{n}
  \DeclareMathSymbol{\Sha}{\mathalpha}{cyrletters}{"58}%
}{}
\makeatother

\newcommand{\Q}{\mathbb{Q}}
\newcommand{\ap}{a_p}
\newcommand{\Omeg}{\Omega_E}
\newcommand{\Tam}{\textstyle\prod_{p \mid N} c_p}
\newcommand{\Tors}{|E(\Q)_{\mathrm{tors}}|}
\newcommand{\Shaan}{|\Sha|_{\mathrm{an}}}
\DeclareMathOperator{\Cov}{Cov}
\DeclareMathOperator{\Corr}{Corr}
\graphicspath{{figures/}}
\renewcommand{\arraystretch}{1.1}

\title[BSD Invariants and Murmurations]{BSD Invariants and Murmurations of Elliptic Curves}

\author{Dane Wachs}
\address{The University of Arizona, Tucson, AZ}
\email{wachs@arizona.edu}
\date{\today}

\begin{document}

\begin{abstract}
We study how the murmurations of elliptic curves over $\Q$ depend on the invariants in the Birch--Swinnerton-Dyer formula, using the Frobenius traces at all primes up to $10^5$ for one curve in each isogeny class of conductor below $10^5$. In these data, at the murmuration scale the rank-$0$ profile depends on the central value $L(E,1)$, in its shape and in the position of its sign changes. Through the series for $L(E,1)$, the covariance of $a_p$ with $L(E,1)$ is exactly a sum of covariances of $a_p$ with the coefficients $a_n$ weighted by a factor that depends on $n$ and the conductor; numerically, in this conductor range, the terms with $n\le50$ reproduce its shape (correlation $0.998$, largest deviation $0.29$ against a largest absolute covariance of $3.03$). Empirically, for $q=3$, $5$, $7$ and $11$ the unweighted covariance of $a_p$ with $a_q$ has nearly the shape of the rank-$0$ profile dilated by a factor near $q$; for $q=2$ the agreement is weaker. Weighting by the Tamagawa product changes the rank-$0$ profile substantially on $0.025<y\le0.125$; after standardisation to the deciles of $\log L(E,1)$, the change summed over these bins is $13\%$ of the unstandardised one ($17\%$ for sums of absolute values), and in the lowest bin the standardised change is larger than the unstandardised one and of opposite sign. This version replaces the earlier versions (v1 and v2); their claim that the Tate--Shafarevich group carries information about Frobenius traces beyond the other invariants taken together cannot hold, because its analytic order is defined by the BSD formula, and it is withdrawn.
\end{abstract}

\maketitle

\section{Introduction}\label{sec:intro}

\subsection{Murmurations}
He, Lee, Oliver and Pozdnyakov averaged the Frobenius trace $\ap$ over elliptic curves over $\Q$ of fixed rank with conductor in a
fixed range: ``Plotting this average as $p$ varies over the primes yields a striking oscillating pattern, the details of which vary with
the rank'' \cite[Abstract]{HLOP}. They note ``a certain invariance to the oscillations as the intervals are scaled to include larger conductors'', which they attribute
to Sutherland \cite[\S1]{HLOP}; plotted against $y=p/N$, the averages behave approximately as functions of $y$. We call
the mean of $\ap$ over rank-$r$ curves, as a function of $y$, the rank-$r$ murmuration profile. For modular forms, as Sawin and
Sutherland put it, ``Convergence of similar averages was proved by Zubrilina \cite{Zub} for modular forms of weight 2 (of which
elliptic curves form a thin subset)'' \cite[Abstract]{SS}; they give a conjecture for elliptic curves ordered by height and, in their
\S1.5, a prediction for curves of prime conductor ordered by conductor, which we compare with our data in Section~\ref{sec:ss}. Martin
\cite{Mar} studies variations, among them averages of $\ap$ weighted by central $L$-values (Section~\ref{sec:prior}).

\subsection{BSD invariants: the definitional point}\label{sec:def}
For a curve $E/\Q$ of conductor $N$ we write $L(E,1)$ for the value at $s=1$ of its $L$-function, $\Omeg$ for the real period, $\Tam$ for
the Tamagawa product, $\Tors$ (or $|T|$) for the torsion order and $\Shaan$ for the analytic order of the Tate--Shafarevich group, all as
stored in Cremona's tables \cite{ecdata}. The first public version of this paper \cite{W2v1} (``v1'') stratified rank-$0$ curves by these invariants
and stated that $|\Sha|$ ``encodes information about the distribution of Frobenius traces at good primes that is not captured by any
other standard BSD invariant''.

At rank $0$ the analytic order is defined by
\begin{equation}\label{eq:shaan}
  \Shaan := \frac{L(E,1)\,\Tors^2}{\Omeg\,\Tam},
\end{equation}
so that
\begin{equation}\label{eq:loglin}
  \log\Shaan = \log L(E,1) + 2\log\Tors - \log\Omeg - \log\Tam .
\end{equation}
Thus $\log\Shaan$ is an exact linear combination of the other four logarithms, $\Shaan$ is a function of $(L(E,1),\Tors,\Omeg,\Tam)$, and
no statistic of the traces can carry information about $\Shaan$ that is not carried by these four quantities taken together. Read with
the other invariants taken jointly, the claim quoted above cannot hold, and it is withdrawn. The sentence of v1 can also be read one
invariant at a time. Under that reading \eqref{eq:shaan} does not exclude it, since a function of four invariants can carry information
that none of them carries alone; but it is then not a statement about $\Sha$, because whatever $\Shaan$ carries beyond each single
invariant it carries as the function \eqref{eq:shaan} of all four. In the tables, \eqref{eq:loglin} holds to $1.0\cdot10^{-14}$ on the $267{,}565$ rank-$0$
curves with $N<10^5$, and $\Shaan$ is an integer square up to rounding on all $1{,}170{,}876$ rank-$0$ curves.

Babei, Banwait, Fong, Huang and Singh make the same point for the LMFDB data. After a logistic regression on the logarithms of the BSD
features classifies perfectly (their item~(2)) and an ordinary least-squares fit detects the linear relation among the logarithms that
their equation~(2.1) implies (item~(3)), they write: ``We do not, however, claim that this is in any way evidence for BSD, since the LMFDB sha\_order field
is computed via Equation (2.1)'' \cite[\S2.1, item (3)]{BBFHS}.

\begin{remark}[isogeny classes]\label{rem:isogeny}
Isogenous curves have the same $L$-function (an isogeny induces an isomorphism of the rational
$\ell$-adic Tate modules as Galois representations, so the local factors agree at every prime, good or bad); hence every coefficient
$a_n$, including $\ap$ at the bad primes, and $L(E,1)$ are invariants of the isogeny class. The factors $\Omeg$, $\Tam$, $\Tors$ and $\Shaan$ are not: in Cremona's tables $88{,}741$ of the $666{,}243$
isogeny classes with more than one curve contain curves whose $\Shaan$, rounded to the nearest integer, differ. Babei et al.\ note the
same for their purpose: ``One must necessarily use more features than merely the $a_p$ values, since these are isogeny invariant, while
[$\Sha(E/\Q)$] is not'' \cite[\S2]{BBFHS}. Consequently a curve-level BSD factor can be associated with $\ap$ only across isogeny
classes, and we take one curve per class as the primary unit of observation (Section~\ref{sec:data}).
\end{remark}

\subsection{This version}
This version replaces the earlier versions (v1 and v2). It keeps the question and Cremona's data, adds new data, namely Frobenius traces at all primes
$p\le10^5$ for one curve in each isogeny class with $N<10^5$, so that $y=p/N$ covers $(0,1]$ for every class, and reports the following.
At small primes (Section~\ref{sec:small}) the traces are positively associated with $L(E,1)$, as the local factor of $L(E,1)$ at $p$
suggests, and the other invariants raise the $R^2$ of a log-linear fit by at most $0.005$ on class representatives ($0.0052$ on all
curves). At the murmuration scale the rank-$0$ profile
depends on $L(E,1)$, in shape and in sign pattern (Section~\ref{sec:L}). Through the series for $L(E,1)$ (Remark~\ref{rem:series}), the
covariance of $\ap$ with $L(E,1)$ is a sum over $n$ of covariances of $\ap$ with $w_na_n$, where the weight $w_n$ depends on the conductor;
for $p\nmid n$ the $n$-th term involves an average of $a_{pn}$ reweighted by $w_n$. In this conductor range the terms with $n\le50$
reproduce the shape of the covariance (correlation $0.998$; largest deviation $0.29$ against a largest absolute covariance of $3.03$;
$n\le100$: $0.9999$ and $0.04$); as $w_n$ decays on the scale $n\asymp\sqrt N$, this is a finding for $N<10^5$ (Section~\ref{sec:reduction}(a)).
The unweighted covariances of $\ap$ with $a_q$ show a further empirical pattern: for $q=3,5,7,11$, $\Cov(\ap,a_q\mid y)$ has correlation
at least $0.95$ with the rank-$0$ profile at $qp/N$ averaged over the pairs of each bin; on fixed bins the correlation is largest at
$s=2.90$, $5.10$, $7.10$, $11.05$, to within $0.05$ in each of two conductor windows, and $s=q$ is within $0.01$ of the maximum for $q=5,7,11$, not for
$q=3$ (on classes with $q\nmid N$: not for $q=3,5$; Section~\ref{sec:reduction}(b)). On rank~$0$, after standardisation to the deciles of $\log L(E,1)$, the
change of the profile under the Tamagawa weight, summed over $0.025<y\le0.125$, is $13\%$ of the unstandardised change ($17\%$ for sums of
absolute values), a ratio of two profile summaries rather than the share of a separately identified effect, and in the lowest bin it
reverses sign ($+0.125\pm0.013$ against $-0.106$); Section~\ref{sec:weights} also treats the real-period weight. Section~\ref{sec:ss} compares prime-conductor data with the prediction of Sawin and Sutherland. Section~\ref{sec:prior} relates the results
to the literature, lists open questions, and summarises the corrections to the earlier versions; the full list is the ancillary file
\texttt{anc/CORRECTIONS\_V1.md}.

\section{Data}\label{sec:data}

\emph{BSD invariants.} Cremona's \texttt{allbsd} tables \cite{ecdata}: $3{,}064{,}705$ curves with $11\le N\le 499{,}998$,
with rank, $\Tors$, $\Tam$, $\Omeg$, $L^{(r)}(E,1)/r!$ and $\Shaan$ as stored. The root number is not in the data; we use
$(-1)^{\mathrm{rank}}$, and PARI's \texttt{ellrootno} gives the same value on all $437{,}226$ class representatives defined below
(\texttt{anc/scripts/s41\_rootno.py}).

\emph{Traces at small primes.} The matrix of $\ap$ at the $500$ primes $p\le 3571$ for all $657{,}396$ curves with $N<10^5$ used in
v1 (at $p\mid N$ the stored entry is the trace of the reduced curve, $0$ or $\pm1$).

\emph{Traces at all $p\le10^5$.} For one curve in each of the $437{,}226$ isogeny classes with $N<10^5$, namely curve number~$1$ in
Cremona's labelling, we computed $\ap$ at all $9{,}592$ primes $p\le10^5$ with PARI/GP~2.17 (\texttt{ellinit} on Cremona's
$a$-invariants, then \texttt{ellap}) \cite{PARI}. At $p\le3571$ these new traces agree with the stored matrix in all
$218{,}613{,}000=437{,}226\times500$ entries (0 mismatches). The bin sums below were checked by brute force: the $K=40$ profiles and the
root-number-weighted mean recomputed directly from the traces with a different exact binning agree to $0$, and for the $11{,}693$ classes
with $N\le3571$ the integer bin sums recomputed from the stored matrix agree exactly (Appendix~\ref{app:comp}).

\emph{Class representatives.} We call these $437{,}226$ curves the \emph{class representatives}: $168{,}760$ of rank $0$, $222{,}439$ of
rank $1$, $45{,}496$ of rank $2$ and $531$ of rank $3$. $\ap$ and $L(E,1)$ are those of the class; $\Omeg$, $\Tam$, $\Tors$ and $\Shaan$ are
those of curve~$1$ (Remark~\ref{rem:isogeny}). We make no claim that curve~$1$ is optimal.

\emph{Window $W$.} For the small-prime analysis: rank $0$, $10^4\le N\le 5\cdot10^4$; $108{,}962$ curves in $68{,}322$ isogeny classes.

\emph{Bins and statistics.} We write $y=p/N$ and use $K$ bins, bin $j=(j/K,(j+1)/K]$, with $K=40$ (width $0.025$) or $K=100$ (width
$0.01$); primes $p>N$ are discarded. Unless stated otherwise only good primes $p\nmid N$ enter. For a set of curves the mean of $\ap$ in a
bin is pairs-weighted, $\sum_E S_{1,E}/\sum_E n_E$, where $S_{1,E}$ is the sum of $\ap$ and $n_E$ the number of good primes of $E$ in the
bin. A class thus enters with weight roughly proportional to $N/\log N$, and large conductors dominate: for the rank-$0$ profile at $K=100$ the
pairs-weighted quartiles of $N$ are $47{,}838$, $69{,}480$ and $85{,}974$ (\texttt{anc/scripts/s46\_pairs\_N.py}). The standard error
(SE) of such a mean is clustered by curve (on class representatives, by isogeny class):
$\mathrm{SE}^2=\sum_E (S_{1,E}-n_E\,\bar a)^2/(\sum_E n_E)^2$. Regression $t$-statistics use errors clustered by isogeny class; on class representatives (one curve per class) this is the
heteroskedasticity-robust HC1 error. Weighted
averages (Section~\ref{sec:weights}) use a bootstrap over classes. We write $z=\text{estimate}/\mathrm{SE}$. The data are
complete finite conductor tables, not random samples, and their averages are exact; the class-bootstrap and cluster standard errors
describe variation under a model in which isogeny classes are exchangeable draws, not uncertainty about the computed finite-table
averages. These standard errors assume
no correlation between isogeny classes, while quadratic twists share $|\ap|$ at good primes. Clustering instead on the $202{,}080$ families of
classes linked by a common $j$-invariant (they contain the quadratic twists; the largest has $618$ classes) changes the standard errors of
Table~\ref{tab:gap40} and Figure~\ref{fig:cov} by less than $2\%$ and leaves the Bonferroni counts of Section~\ref{sec:L} unchanged
(\texttt{anc/csv/twist\_cluster.json}).

\section{Small primes}\label{sec:small}

In this section the sample is the window $W$, and at each prime $p$ the curves with $p\mid N$ are dropped.

An ordinary least-squares fit of $\ap$ on $[1,\,N,\,\log L(E,1)]$ has $R^2=0.240$, $0.250$, $0.172$, $0.122$ at
$p=2,3,5,7$ on class representatives (Table~\ref{tab:r2}). Adding $\log\Omeg$, $\log\Tam$, $\log\Tors$ and the sign of the minimal
discriminant raises $R^2$ by at most $0.005$ at every $p\le23$ on class representatives, and by at most $0.0052$ on all curves. By
\eqref{eq:loglin}, any four of the five logarithms $\log L(E,1)$, $\log\Omeg$, $\log\Tam$, $\log\Tors$, $\log\Shaan$ span the same space;
the coefficients of $\log\Omeg$, $\log\Tam$ and $\log\Tors$ depend on which four of the five logarithms are used, and we do not interpret
them. In the fit on $[1,\,N,\,\log L(E,1)]$ alone the coefficient of $\log L(E,1)$ is positive at every $p\le23$, from $+0.77$ ($p=2$)
to $+1.24$ ($p=7$) on class representatives, with $t=40$ to $91$ (standard errors clustered by isogeny class, which on class
representatives is HC1; \texttt{anc/c4\_sign.json}).

Heuristically, the positive association of $\ap$ with $L(E,1)$ at small $p$ is what the Euler product suggests: formally
$L(E,1)=\prod_{p\nmid N}(1-\ap p^{-1}+p^{-1})^{-1}\prod_{p\mid N}(1-\ap p^{-1})^{-1}$, a product that does not converge absolutely at
$s=1$, in which a larger $\ap$ at a small good prime enlarges the factor at $p$. This heuristic is stated by Knightly and Reno for weights
given by central $L$-values: ``The weight $w_f$ depends directly on $\lambda_f(p)$ via the local $L$-factor [\dots]. Assuming the
remaining $L$-factors do not affect the distribution of the $\lambda_f(p)$, [\dots]'' \cite[Remark~3.4]{KR}. For weight $k>2$ they prove
the corresponding statement \cite[Theorem~1.5]{KR} ($\chi$ primitive and real of conductor $D$ coprime to $N$, $p\nmid DN$ fixed, $N>1$):
weighted by $\Lambda(\frac12,f\times\chi)/\|f\|^2$, the normalised Hecke eigenvalues
$\lambda_f(p)$ of an orthogonal basis of $T_p$-eigenforms of $S_k(N)$ ``become equidistributed in $[-2,2]$ with respect to the probability measure''
$p\,((p+1)-x\chi(p)\sqrt p)^{-1}\,d\mu_\infty(x)$ as $N+k\to\infty$, where $\mu_\infty$ is the Sato--Tate measure; for trivial $\chi$ and
$x=\ap/\sqrt p$ the density is the factor $(1-\ap p^{-1}+p^{-1})^{-1}$ above. Earlier, Ramakrishnan and Rogawski stated, for even
$k>2$, prime level $N$ and a quadratic character $\chi$ of conductor $D<0$ with $\chi(-N)=1$, an asymptotic formula for the sum of
$L(\frac12,\varphi\otimes\chi)L(\frac12,\varphi)/(\varphi,\varphi)$ over the newforms $\varphi$ whose normalised $a_p(\varphi)$ lies in a
given interval, $p\nmid ND$ fixed, with main term proportional to an explicit measure $\mu_p$ of the interval and error
$O(N^{-k/2+\varepsilon})$ \cite[Theorem~A]{RR} (see also their erratum \cite{RRe}, whose text we could not access); Knightly and Reno
remark that the argument for characteristic functions of intervals with this error term is incomplete, because the error for an
approximating polynomial depends on the polynomial \cite[Remark~5.3(1)]{KR}. Ramakrishnan and Rogawski
write ``For simplicity of exposition, we are not treating the case $k=2$ here'' \cite[\S1]{RR}. Our setting differs (weight $2$, elliptic
curves, conditioning on $L(E,1)$ rather than weighting, finite $N$), so these results do not apply to it directly.

\begin{table}[ht]
\centering
\caption{$R^2$ of per-prime ordinary least-squares fits (one fit per prime, no binning in $y$) in the window $W$ (rank $0$, $10^4\le N\le5\cdot10^4$); at each $p$ the
curves with $p\mid N$ are dropped. $R^2_L$: $\ap$ on $[1,N,\log L(E,1)]$. $R^2_{\mathrm{full}}$: in addition $\log\Omeg$, $\log\Tam$,
$\log\Tors$ and $\mathbf 1(\Delta_{\min}>0)$. $\Delta R^2=R^2_{\mathrm{full}}-R^2_L$. Left: class representatives (curve number~$1$);
right: all curves. No standard errors enter this table. Source: \texttt{csv/contract\_numbers.json} (``A''), script
\texttt{contract\_numbers.py}.}\label{tab:r2}
{\small
\begin{tabular}{r ccc ccc}
\toprule
 & \multicolumn{3}{c}{class representatives ($68{,}322$)} & \multicolumn{3}{c}{all curves ($108{,}962$)}\\
\cmidrule(lr){2-4}\cmidrule(lr){5-7}
$p$ & $R^2_{L}$ & $R^2_{\mathrm{full}}$ & $\Delta R^2$ & $R^2_{L}$ & $R^2_{\mathrm{full}}$ & $\Delta R^2$\\
\midrule
2 & 0.2398 & 0.2431 & 0.0033 & 0.2151 & 0.2196 & 0.0045\\
3 & 0.2501 & 0.2551 & 0.0050 & 0.2201 & 0.2253 & 0.0052\\
5 & 0.1724 & 0.1769 & 0.0044 & 0.1734 & 0.1780 & 0.0047\\
7 & 0.1217 & 0.1250 & 0.0033 & 0.1269 & 0.1305 & 0.0035\\
11 & 0.0674 & 0.0701 & 0.0027 & 0.0746 & 0.0773 & 0.0027\\
13 & 0.0584 & 0.0596 & 0.0012 & 0.0624 & 0.0638 & 0.0014\\
17 & 0.0400 & 0.0417 & 0.0017 & 0.0428 & 0.0447 & 0.0019\\
19 & 0.0311 & 0.0325 & 0.0014 & 0.0345 & 0.0361 & 0.0016\\
23 & 0.0244 & 0.0260 & 0.0016 & 0.0272 & 0.0286 & 0.0014\\
\bottomrule
\end{tabular}
}
\end{table}

\section{The murmuration scale}\label{sec:scale}

\subsection{The profiles}\label{sec:profiles}
\begin{figure}[ht]
\centering
\includegraphics[width=0.88\textwidth]{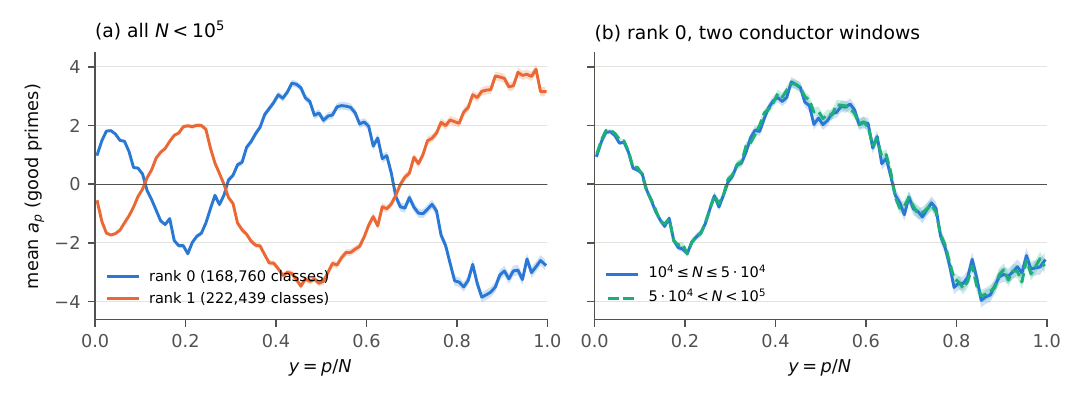}
\caption{Murmuration profiles at $K=100$: mean of $\ap$ over good primes in each bin of $y=p/N$, class representatives, bands $\pm2$
curve-cluster SE. (a) Rank $0$ and rank $1$, all $N<10^5$. (b) Rank $0$ on $10^4\le N\le5\cdot10^4$ and on $5\cdot10^4<N<10^5$.
Source: \texttt{data\_pN/results/pnA\_profiles\_K100.csv}.}\label{fig:profiles}
\end{figure}

Figure~\ref{fig:profiles} shows the rank-$0$ ($168{,}760$ classes) and rank-$1$ ($222{,}439$ classes) profiles on
$y\in(0,1]$ at $K=100$. The rank-$0$ profile changes sign at $y\approx0.111$, $0.292$ and $0.660$ (linear interpolation between bin
centres); its maximum is $+3.44$ near $y=0.44$ and its minimum $-3.85$ near $y=0.86$. Computed separately on $N\in[10^4,5\cdot10^4]$ and
on $N\in(5\cdot10^4,10^5)$, the two rank-$0$ profiles have correlation $0.998$ over the $100$ bins.

\subsection{Conditioning on \texorpdfstring{$L(E,1)$}{L(E,1)}}\label{sec:L}
We split the rank-$0$ class representatives into quartiles $Q_1,\dots,Q_4$ of $\log L(E,1)$, with $Q_1$ the smallest values, and
compare the profiles of $Q_4$ and $Q_1$ bin by bin. Unless stated otherwise the comparison is stratified by conductor: the quartiles are
formed within each of $20$ conductor bins, so that each quartile is approximately balanced across the conductor bins. The conductor bins have as edges the quantiles of $\log N$ at $k/20$
($k=1,\dots,19$; linear interpolation, as in \texttt{numpy.quantile}) over the rank-$0$ class representatives, and a class whose $\log N$
equals an edge goes to the lower bin; within a conductor bin the quartile edges are the quantiles at $1/4$, $1/2$, $3/4$, and a value equal
to an edge goes to the upper quartile. With these rules Table~\ref{tab:gap40} is reproduced from the bin sums by
\texttt{anc/scripts/s23\_twist\_cluster.py}. For the difference $Q_4-Q_1$ in a bin, $z$
is the difference divided by the two groups' curve-cluster SEs combined in quadrature.

\begin{figure}[ht]
\centering
\includegraphics[width=0.88\textwidth]{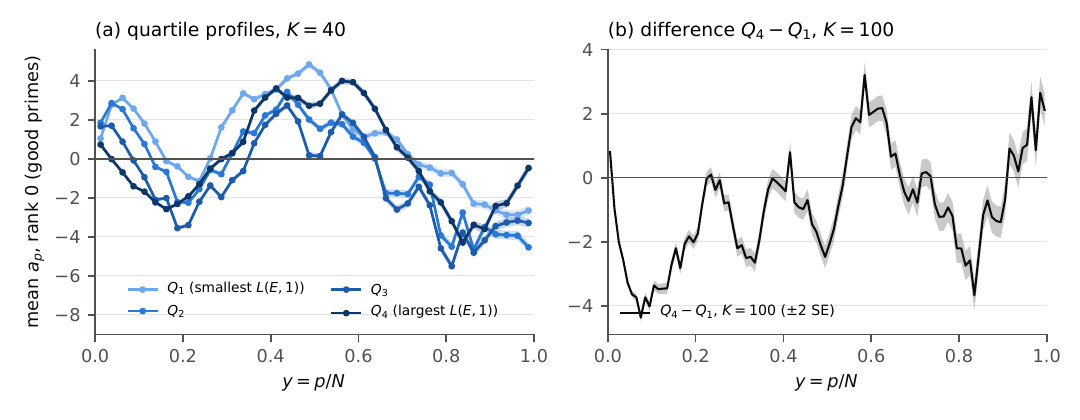}
\caption{Rank-$0$ class representatives by quartile of $\log L(E,1)$, quartiles formed over all $N$ ($42{,}190$ classes each); good
primes. (a) Mean $\ap$ per quartile at $K=40$, bands $\pm2$ curve-cluster SE. (b) $Q_4-Q_1$ at $K=100$, band $\pm2$ SE (the two groups'
curve-cluster SEs in quadrature). Table~\ref{tab:gap40} gives $Q_4-Q_1$ with quartiles formed within conductor bins. Sources:
\texttt{data\_pN/results/pnC\_conditioned\_K40.csv}, \texttt{pnC\_conditioned\_K100.csv}.}\label{fig:quartiles}
\end{figure}

At $K=40$ (Table~\ref{tab:gap40}):
\begin{itemize}[leftmargin=1.5em]
\item $|z|\ge3.23$ (the Bonferroni threshold for $40$ bins at two-sided level $5\%$) at $31$ of the $40$ bins.
\item On the Bonferroni-significant bins, $Q_4-Q_1$ is negative for $0<y\le0.225$ (largest gap $-4.01$ at $y\approx0.09$, $z=-97$),
negative again on $0.25$--$0.375$ and $0.425$--$0.525$, positive on $0.55$--$0.65$ ($+2.46$ at $y\approx0.59$, $z=+20$) and on
$0.925$--$1$, and negative on $0.75$--$0.90$.
\item At $K=100$ the lowest bin, $y\le0.01$, is positive: $+0.80$, $z=+73$ (Figure~\ref{fig:quartiles}(b)). Here the quartiles are
formed over all $N$, not within conductor bins; the two definitions differ by at most $0.061$ at $y\le0.125$ (Section~\ref{sec:robust}).
This gap is carried by primes $p=O(\sqrt N)$ and is not evidence about the murmuration scale (see the end of this subsection).
\end{itemize}
The quartile profiles themselves (Figure~\ref{fig:quartiles}(a); quartiles over all $N$) differ in shape and in the positions of their
sign changes, not only in amplitude: the lowest zero crossing lies near $y\approx0.16$ for $Q_1$ and near $y\approx0.04$ for $Q_4$.
Without the good primes $p\le2\sqrt N$ these crossings stay at $0.1595$ ($Q_1$) and move from $0.0371$ to $0.0370$ ($Q_4$, whose lowest
$K=40$ bin stays positive, $+0.52$, $z=42$; \texttt{anc/scripts/s60\_lowbin\_split.py}).

\begin{figure}[ht]
\centering
\includegraphics[width=0.85\textwidth]{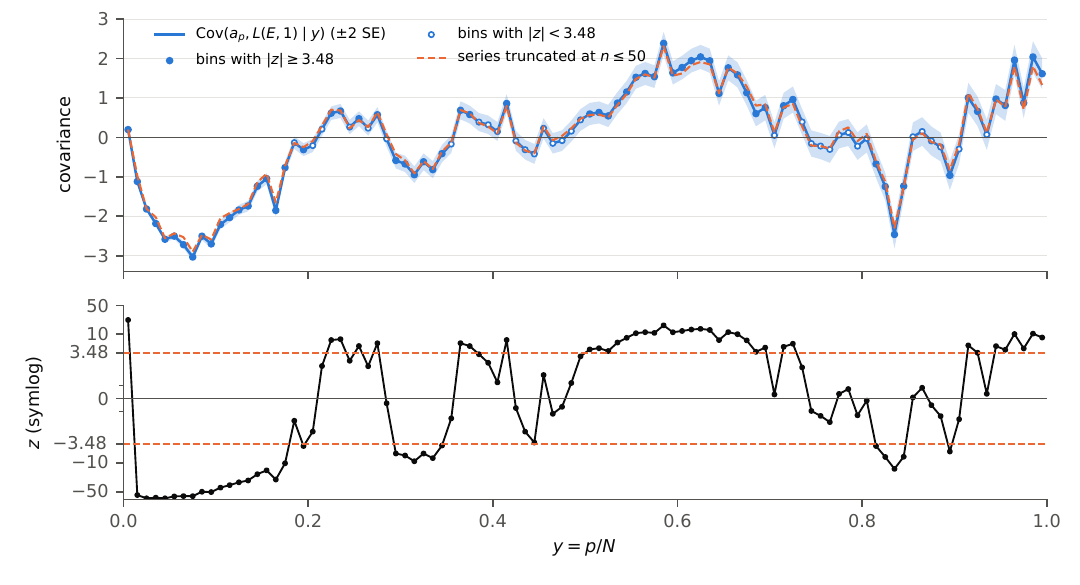}
\caption{$\Cov(\ap,L(E,1)\mid y)$ at $K=100$: rank-$0$ class representatives, good primes, pairs-weighted; band $\pm2$ linearised
curve-cluster SE. Filled markers: $|z|\ge3.48$ (Bonferroni, $100$ bins, two-sided $5\%$); open markers: $|z|<3.48$. Dashed: the covariance
with the series \eqref{eq:series} truncated at $n\le50$ (Section~\ref{sec:reduction}). Lower panel: $z$ on a symmetric-log scale, dashed
lines at $\pm3.48$. Source: \texttt{anc/csv/cov\_L\_K100.csv} (script \texttt{s20\_afe\_cov.py}); figure \texttt{f03\_cov\_L.py}.}\label{fig:cov}
\end{figure}

The same dependence appears in the covariance. Let $\Cov(\ap,L(E,1)\mid y)$ be the pairs-weighted covariance over the
(curve, good prime) pairs of a bin, rank-$0$ class representatives, $K=100$, with a linearised curve-cluster SE (Figure~\ref{fig:cov};
\texttt{anc/csv/cov\_L\_K100.csv}). It has $|z|\ge3.48$ (Bonferroni over $100$ bins) at $67$ of the $100$ bins. It is positive at
$y\le0.01$ ($+0.20$, $z=+22$; this bin is treated at the end of this subsection) and negative on all of $0.01<y\le0.21$, where it changes
sign ($y\approx0.210$ by linear interpolation);
on $0.01<y\le0.18$ every bin is Bonferroni-significant (minimum $-3.03$ at $y\approx0.075$, $z=-63$), and of the next three bins only
$0.19<y\le0.20$ is (\texttt{anc/scripts/s47\_covL\_negrun.py}). Further, on runs of
Bonferroni-significant bins, it is negative on $0.29$--$0.35$ and $0.81$--$0.85$ and positive on $0.50$--$0.70$ and $0.94$--$1$; twelve
further significant bins lie in runs of at most two bins. The correlations are small, $|\Corr|\le0.023$ in every bin; the $z$-scores are
large because the sample is large (about $1.5\cdot10^7$ pairs in the lowest bin). With $\log L(E,1)$ in place of $L(E,1)$ the covariance has
correlation $0.89$ with this one over the $100$ bins and $72$ Bonferroni-significant bins (same file). At $K=100$ the sign change nearest
the origin lies between the first two bin centres ($y\approx0.007$ by linear interpolation); it rests on primes $p=O(\sqrt N)$ (below).
With bins of width $0.001$ on $y\le0.03$,
computed from the traces (\texttt{anc/scripts/s33\_crossing\_drift.py}), the first sign change (linear interpolation between bin centres)
lies at $y=0.0068$, $0.0060$, $0.0055$, $0.0053$ and $0.0045$ in the conductor windows $N<2\cdot10^4$, $[2,4)\cdot10^4$, $[4,6)\cdot10^4$,
$[6,8)\cdot10^4$ and $[8,10)\cdot10^4$ (between $33{,}285$ and $34{,}309$ classes each), with $95\%$ class-bootstrap intervals
$[0.0062,0.0078]$, $[0.0057,0.0064]$, $[0.0044,0.0059]$, $[0.0044,0.0057]$ and $[0.0041,0.0050]$ ($B=1000$); over all $N<10^5$ it lies at
$0.0052$ $[0.0048,0.0056]$. The crossing is a linearly interpolated zero of a curve that is nearly flat near it in some windows, so its
bootstrap distribution is lumpy and these percentile intervals are rough. The crossing decreases from the lowest to the highest window;
the intervals of adjacent windows overlap, those of the two extreme windows do not, so for $N<10^5$ the first sign change does not sit
at a fixed $y$.

\emph{The lowest bins: primes $p=O(\sqrt N)$.} The bin $y\le0.01$ at $K=100$ holds the primes $p\le N/100$. For $N\ge10^4$ these include
every $p\le\sqrt N$, and at the pairs-weighted conductor quartiles of Section~\ref{sec:data} the upper edge $N/100$ lies at $2.19\sqrt N$,
$2.64\sqrt N$ and $2.93\sqrt N$. The positive lowest-bin covariance, the lowest-bin quartile gap and the $K=100$ sign change nearest
the origin are carried by primes $p=O(\sqrt N)$, the regime of the Euler-factor heuristic of Section~\ref{sec:small}, and are not
evidence about the murmuration scale $p\asymp N$ (\texttt{anc/scripts/s60\_lowbin\_split.py}, which reproduces the bin sums, the
covariance and the crossings above).
\begin{itemize}[leftmargin=1.5em]
\item \emph{Covariance.} For a good prime $p$ the relations in Remark~\ref{rem:series} give $\ap L(E,1)=\sum_nw_na_{pn}+h(E,p)$, where
the Hecke part $h(E,p)=\sum_{m\ge1}p\,w_{pm}a_m=\sum_{m\ge1}2m^{-1}a_me^{-2\pi pm/\sqrt N}$ collects the terms $p\,a_{n/p}$ at $p\mid n$.
Its term $m=1$, $2e^{-2\pi p/\sqrt N}$, comes from $a_p^2=a_{p^2}+p$ and is a deterministic function of $p$ and $N$ (the self-term), and
$h$ is non-negligible only for $p=O(\sqrt N)$. In the bin $y\le0.01$ the Hecke part contributes $0.184$ to the covariance (the $0.184$
of Section~\ref{sec:reduction}), $0.172$ of it from the self-term; the remaining family part,
$\Cov_b(\ap,L(E,1))-M_b(h)=M_b(\sum_nw_na_{pn})-M_b(\ap)M_b(L(E,1))$, is $+0.015\pm0.009$ ($z=1.6$). In every other bin the Hecke part
is below $1.5\cdot10^{-4}$. Without it $66$ of the $100$ bins are Bonferroni-significant (only the lowest bin changes), and the family part
has no significant positive bin below $y=0.22$. So the sign change nearest the origin at $K=100$ rests on the Hecke part: without it the
lowest bin is still positive, but not significantly.
\item \emph{Quartile gap.} Without the good primes $p\le\sqrt N$ ($48\%$ of the pairs of the bin) the $K=100$ gap $Q_4-Q_1$ in the bin
$y\le0.01$ falls from $+0.80$ ($z=73$) to $+0.24\pm0.02$ ($z=11.8$); without $p\le2\sqrt N$ ($82\%$ of the pairs) it is $-0.04\pm0.04$
($z=-1.0$). Removing these primes empties the bin for every class with $N<10{,}101$ (respectively $N<40{,}104$). On exactly the classes
that keep a prime in the bin, the gap with all primes is $+0.79$ ($z=72$) and $+0.76$ ($z=59$), so the fall comes from the primes
removed, not from the change of classes. Since for $N<10^5$ the bin contains no prime above $3.2\sqrt N$, these cuts in $p$ are also cuts in $y$;
they locate the gap at small $p$ but do not decide how it scales (compare open question~2).
\item \emph{Unaffected.} Beyond the lowest bin the Hecke part is below $1.5\cdot10^{-4}$, so the negative run and the sign pattern of
the covariance stand. Without the good primes $p\le2\sqrt N$ the gaps of Table~\ref{tab:gap40} on $y>0.025$ change by at most $0.0013$
and their $z$ by at most $0.16$, with no change of sign or of Bonferroni significance (the classes with such primes there have
$N\le6{,}276$); Figure~\ref{fig:quartiles}(a) is as stated above. The drift of the first crossing holds for the family part. In the bins
of width $0.001$ the family part is positive on $0.001<y\le0.004$ in every window ($z\ge4.0$), and its first sign change lies at
$0.0060$, $0.0055$, $0.0053$ and $0.0045$ in the four upper windows ($95\%$ intervals $[0.0056,0.0064]$, $[0.0044,0.0059]$,
$[0.0044,0.0057]$, $[0.0041,0.0050]$, same bootstrap) and at $0.0051$ $[0.0047,0.0056]$ over all $N<10^5$. In the lowest window the
family part is $-0.006$ ($z=-0.4$) in the first fine bin (there $p\le19$), so the rule above puts its first sign change between the
first two fine bins; the first change to a run of at least five negative bins, a secondary rule fixed in advance in
\texttt{s33\_crossing\_drift.py}, lies at $0.0065$ $[0.0061,0.0097]$, and in the other windows the two rules agree. The intervals of
the extreme windows remain disjoint.
\item \emph{Scale.} The positive region below the first crossing lies at $p$ of order $\sqrt N$. At each window's pairs-weighted median
conductor $N_{\mathrm{med}}=13{,}790$, $31{,}338$, $50{,}850$, $70{,}650$, $90{,}480$ (weights: good primes in the fine bins), the family
part's first crossing $y^*$ corresponds to $p^*/\sqrt{N_{\mathrm{med}}}=y^*\sqrt{N_{\mathrm{med}}}=0.76$, $1.06$, $1.23$, $1.41$,
$1.36$ (lowest window: secondary rule; $0.80$, $1.07$, $1.24$, $1.41$, $1.37$ for the full covariance). From the lowest to the highest
window $N_{\mathrm{med}}$ grows by a factor $6.6$, $y^*$ falls by a factor $1.4$ and $p^*/\sqrt{N_{\mathrm{med}}}$ rises by a factor
$1.8$; the least-squares slope of $\log y^*$ on $\log N_{\mathrm{med}}$ is $-0.17$ (a fixed $y$ gives $0$, a fixed $p/\sqrt N$ gives
$-1/2$). Over these windows neither $p/N$ nor $p/\sqrt N$ is fixed (the $95\%$ intervals of $y^*\sqrt{N_{\mathrm{med}}}$ in the lowest and the highest window, $[0.72,1.14]$ and $[1.24,1.49]$, are disjoint, as are those of $y^*$).
\end{itemize}

\begin{table}[ht]
\centering
\caption{$Q_4-Q_1$ at $K=40$: difference of the rank-$0$ mean $\ap$ (good primes) between the largest and the smallest quartile of
$\log L(E,1)$, quartiles formed within $20$ conductor bins; class representatives ($168{,}760$ rank-$0$ classes). Bin $j$ is
$y\in(j/40,(j+1)/40]$; the centre is shown. $z$: the difference over the two groups' curve-cluster SEs combined in quadrature.
${}^{*}$: $|z|\ge3.23$ (Bonferroni, $40$ bins, two-sided $5\%$). Source: \texttt{csv/contract\_numbers.json} (``B''), script
\texttt{contract\_numbers.py}. $|z|\ge10$ is rounded to an integer.}\label{tab:gap40}
{\small
\begin{tabular}{rcrr@{\hspace{2.2em}}rcrr}
\toprule
bin & $y$ (centre) & $Q_4-Q_1$ & $z$ & bin & $y$ (centre) & $Q_4-Q_1$ & $z$\\
\midrule
0 & 0.0125 & $-0.335$ & $-31$$^{*}$ & 20 & 0.5125 & $-1.473$ & $-13$$^{*}$\\
1 & 0.0375 & $-2.839$ & $-117$$^{*}$ & 21 & 0.5375 & $+0.067$ & $+0.6$\\
2 & 0.0625 & $-3.856$ & $-115$$^{*}$ & 22 & 0.5625 & $+1.791$ & $+15$$^{*}$\\
3 & 0.0875 & $-4.005$ & $-97$$^{*}$ & 23 & 0.5875 & $+2.457$ & $+20$$^{*}$\\
4 & 0.1125 & $-3.500$ & $-74$$^{*}$ & 24 & 0.6125 & $+2.285$ & $+18$$^{*}$\\
5 & 0.1375 & $-3.085$ & $-59$$^{*}$ & 25 & 0.6375 & $+1.369$ & $+11$$^{*}$\\
6 & 0.1625 & $-2.433$ & $-41$$^{*}$ & 26 & 0.6625 & $+0.283$ & $+2.2$\\
7 & 0.1875 & $-1.849$ & $-29$$^{*}$ & 27 & 0.6875 & $-0.235$ & $-1.8$\\
8 & 0.2125 & $-0.906$ & $-13$$^{*}$ & 28 & 0.7125 & $-0.073$ & $-0.5$\\
9 & 0.2375 & $-0.106$ & $-1.5$ & 29 & 0.7375 & $-0.294$ & $-2.2$\\
10 & 0.2625 & $-0.481$ & $-6.2$$^{*}$ & 30 & 0.7625 & $-1.146$ & $-8.2$$^{*}$\\
11 & 0.2875 & $-1.545$ & $-19$$^{*}$ & 31 & 0.7875 & $-1.360$ & $-9.6$$^{*}$\\
12 & 0.3125 & $-2.100$ & $-25$$^{*}$ & 32 & 0.8125 & $-2.339$ & $-16$$^{*}$\\
13 & 0.3375 & $-2.404$ & $-27$$^{*}$ & 33 & 0.8375 & $-2.910$ & $-20$$^{*}$\\
14 & 0.3625 & $-0.504$ & $-5.4$$^{*}$ & 34 & 0.8625 & $-1.045$ & $-6.9$$^{*}$\\
15 & 0.3875 & $-0.111$ & $-1.2$ & 35 & 0.8875 & $-1.207$ & $-7.8$$^{*}$\\
16 & 0.4125 & $+0.120$ & $+1.2$ & 36 & 0.9125 & $+0.226$ & $+1.4$\\
17 & 0.4375 & $-0.935$ & $-9.0$$^{*}$ & 37 & 0.9375 & $+0.593$ & $+3.7$$^{*}$\\
18 & 0.4625 & $-1.180$ & $-11$$^{*}$ & 38 & 0.9625 & $+1.503$ & $+9.2$$^{*}$\\
19 & 0.4875 & $-1.926$ & $-17$$^{*}$ & 39 & 0.9875 & $+2.304$ & $+14$$^{*}$\\
\bottomrule
\end{tabular}
}
\end{table}

\subsection{The covariance with \texorpdfstring{$L(E,1)$}{L(E,1)} rewritten through the series for \texorpdfstring{$L(E,1)$}{L(E,1)}}\label{sec:reduction}
For a finite set of curves and a bin $b$ of $y$, write $M_b(f)$ for the pairs-weighted mean of a function $f(E,p)$ over the pairs $(E,p)$
with $p$ a good prime of $E$ in the bin, and $\Cov_b(\ap,x)=M_b(\ap x)-M_b(\ap)M_b(x)$ for a function $x(E)$; this is the covariance of
Section~\ref{sec:L}.

\begin{remark}\label{rem:series}
Let $\mathcal E$ be a finite set of elliptic curves $E/\Q$ with $L(E,1)\neq0$, let $b$ be a bin of $y$, and put
$w_n(E)=2n^{-1}e^{-2\pi n/\sqrt{N_E}}$, where $N_E$ is the conductor of $E$. Then for every $E\in\mathcal E$
\begin{equation}\label{eq:series}
  L(E,1)=\sum_{n\ge1}w_n(E)\,a_n(E),
\end{equation}
the series converging absolutely, and
\begin{equation}\label{eq:covsum}
  \Cov_b\bigl(\ap,L(E,1)\bigr)=\sum_{n\ge1}\Cov_b(\ap,w_na_n),\qquad \Cov_b(\ap,w_na_n)=M_b(w_n\,\ap a_n)-M_b(\ap)\,M_b(w_na_n).
\end{equation}
For a good prime $p$, $\ap a_n=a_{pn}$ if $p\nmid n$, and $\ap a_n=a_{pn}+p\,a_{n/p}$ if $p\mid n$. If no prime of the bin divides $n$, the
first piece of the $n$-th term of \eqref{eq:covsum}, $M_b(w_n\ap a_n)=M_b(w_na_{pn})$, is an average of $a_{pn}$ over the pairs of the bin
(so $pn/N$ lies in the bin $nb$, at position $ny$), reweighted by $w_n(E)$; this weight depends on $N_E$, so the average is not the
unweighted mean of the profiles. The second piece, $M_b(\ap)M_b(w_na_n)$, is the profile $M_b(\ap)$ at $y$ times the bin-dependent factor
$M_b(w_na_n)$.

This is linearity of a finite covariance applied to Cremona's absolutely convergent series for $L(E,1)$. By modularity $L(E,s)=L(f,s)$
for a newform $f$ of weight $2$ and level $N=N_E$ with Fourier coefficients $a(n,f)=a_n(E)$. Cremona shows that
$L(f,1)=(\varepsilon_N-1)I_f(i/\sqrt N)$, where $\varepsilon_N$ is the eigenvalue of the Fricke involution $W_N$ on $f$, and concludes:
``Thus if $L(f,1)\neq0$, then necessarily $\varepsilon_N=-1$'' \cite[\S2.11, (2.11.1)]{Cre}. His Proposition~2.11.1 states: if
$f\in S_2(N)$ and $f|W_N=-f$ then $L(f,1)=2\sum_{n=1}^\infty\frac{a(n,f)}{n}\exp(-2\pi n/\sqrt N)$ \cite[(2.11.2)]{Cre}; this is
\eqref{eq:series}. The coefficients are multiplicative, $a_{mn}=a_ma_n$ for coprime $m$ and $n$, and the local factors
$(1-\ap p^{-s}+p^{1-2s})^{-1}$ for $p\nmid N$ and $(1-\ap p^{-s})^{-1}$ for $p\mid N$ give $\ap a_{p^k}=a_{p^{k+1}}+p\,a_{p^{k-1}}$ ($k\ge1$)
for good $p$ and $a_{p^k}=\ap^k$ for bad $p$; writing $n=p^km$ with $p\nmid m$ gives the statement about $\ap a_n$. With the Hasse bound
$|\ap|\le2\sqrt p$ the same relations give $|a_n|\le d(n)\sqrt n$, so \eqref{eq:series} converges absolutely. Each $M_b$ is a finite
weighted sum over the pairs of the bin, so it commutes with the sum over $n$, and \eqref{eq:covsum} follows because $\Cov_b(\ap,\cdot)$ is
linear.
\end{remark}

All $1{,}170{,}876$ rank-$0$ curves of the tables have $L(E,1)>0$ (\texttt{anc/csv/l1\_positive.json}), so Remark~\ref{rem:series} applies to the covariance of Section~\ref{sec:L} for
every bin and with no limit in $N$. It is an exact form of the remark of Babei et al.\ that ``the $a_p$ values determine the $L$-function of
the curve and thus already encode the special value'' \cite[\S2.3]{BBFHS}; Cowan studies how murmurations arise from averaged approximate
functional equations \cite{Cow} (Section~\ref{sec:lit}).

Summing \eqref{eq:covsum} over $n$ gives the exact rewriting
\begin{equation}\label{eq:twoparts}
  \Cov_b(\ap,L)=M_b(\ap L)-M_b(L)\,M_b(\ap),\qquad M_b(\ap L)=\sum_{n\ge1}M_b(w_n\ap a_n),\qquad L=L(E,1),
\end{equation}
where $M_b(\ap L)/M_b(L)$ is the $L(E,1)$-weighted mean of $\ap$ (Section~\ref{sec:weights}): a weighted sum of the averages
$M_b(w_n\ap a_n)$ minus a multiple of the profile at $y$ itself. The term $n=1$, $M_b(w_1\ap)$, is not at a composite index; with
$M_b(w_1)M_b(\ap)$ it forms $C_b=\Cov_b(\ap,w_1)$, so that
\begin{equation}\label{eq:split}
  \Cov_b(\ap,L)=A_b-M_b(L-w_1)\,M_b(\ap)+C_b,\qquad A_b=\sum_{n\ge2}M_b(w_n\ap a_n)=M_b\bigl(\ap(L-w_1)\bigr).
\end{equation}
On the rank-$0$ class representatives at $K=100$ (\texttt{anc/scripts/s30\_decomp\_sizes.py}), $M_b(L(E,1)-w_1)$ lies
between $0.8707$ and $0.8731$ in every bin, $\max_b|C_b|=0.0023$, and the root mean square over the $100$ bins is $2.33$ for $A_b$, $1.91$ for
$M_b(L(E,1)-w_1)M_b(\ap)$ and $1.23$ for the covariance, whose correlation with $A_b$ over the bins is $0.57$. The near-constancy of
$M_b(L(E,1)-w_1)$ is a property of the pair weights, not of the traces: $M_b(L(E,1))$ lies between $2.817$ and $2.820$ and $M_b(w_1)$
between $1.946$ and $1.947$ over the $100$ bins, because the conductor mix of the pairs is nearly the same in every bin (the
pairs-weighted median of $N$ lies between $68{,}952$ and $69{,}758$; \texttt{s46\_pairs\_N.py}). The split at $n=1$ is a
choice: in \eqref{eq:twoparts}, where both parts contain the term $n=1$, the root mean squares are $6.36$ and $6.17$. In $A_b$ the terms
with $p\mid n$ are not averages of $a_{pn}$ alone, since $\ap a_n=a_{pn}+p\,a_{n/p}$ there; the extra parts $p\,a_{n/p}$ contribute
$0.184$ to $A_b$ in the bin $y\le0.01$ and less than $1.5\cdot10^{-4}$ in every other bin (\texttt{anc/scripts/s36\_pdivn.py}); this
is the Hecke part of the lowest bins in Section~\ref{sec:L}.
Remark~\ref{rem:series} holds for every finite set of curves and every bin; it does not by itself give the size, the sign pattern or the
zeros of $\Cov(\ap,L(E,1)\mid y)$, which come from the averages $M_b(w_n\ap a_n)$. What the data show about these averages is in (a)
and (b). In the rest of this subsection the sample is the $168{,}760$ rank-$0$ class representatives and $K=100$
(\texttt{anc/scripts/s20\_afe\_cov.py}).

(a) \emph{The series.} With $a_n$ computed from the traces ($a_1=1$; $a_{q^{k+1}}=a_qa_{q^k}-q\,a_{q^{k-1}}$ for $q\nmid N$,
$a_{q^k}=a_q^k$ for $q\mid N$; $a_{mn}=a_ma_n$ for coprime $m$, $n$), the series \eqref{eq:series} truncated at $n\le1000$ reproduces
Cremona's $L(E,1)$ on all $168{,}760$ rank-$0$ class representatives with relative error at most $2.7\cdot10^{-8}$ (median
$2\cdot10^{-14}$). Truncating \eqref{eq:covsum} at $n\le10$, $20$, $50$ and $100$ gives $\Cov(\ap,L(E,1)\mid y)$ with correlation $0.886$,
$0.975$, $0.998$ and $0.9999$ with the full covariance over the $100$ bins, and maximum absolute difference $2.04$, $1.03$, $0.29$ and $0.04$,
against $\max|\Cov|=3.03$ (Figure~\ref{fig:cov}, dashed: $n\le50$). In this conductor range the covariance with $L(E,1)$ at $p\asymp N$
is thus carried by the covariances of $\ap$ with $w_na_n$ for small $n$; since $w_n(E)$ decays on the scale $n\asymp\sqrt{N_E}$, this is
a finding for $N<10^5$, not a statement uniform in $N$. The series identity connects the central-value covariance to these weighted
coefficient covariances, composite $n$ included; it does not reduce it to the unweighted $\Cov(\ap,a_q\mid y)$ of (b), since
$\Cov(X,WY)$ is not determined by $\Cov(X,Y)$. The unweighted prime-index covariances show the further empirical dilation pattern of
(b); the quartile gaps are not covered by the identity (see (c)).

(b) \emph{The composite indices $qp$.} Let
$\Cov(\ap,a_q\mid y)$ be the pairs-weighted covariance of $\ap$ (good primes $p$) with $a_q$, with the stored value of $a_q$ when $q\mid N$
and the linearised curve-cluster SE of Section~\ref{sec:L}; the bins include the pair $p=q$, where $a_qa_p=a_q^2\neq a_{q^2}$, and
removing it changes the correlations at $s=q$ below by at most $0.0014$ and moves no maximiser (variant \texttt{p\_ne\_q}). Let $P_0$ be the rank-$0$ profile and
$P_0(x)$ its value at $x$ by linear interpolation between bin centres. A bin of width $0.01$ in $y$ corresponds to a stretch of width
$0.01s$ of the argument of $P_0(sy)$ ($0.11$ for $s=11$), so we compare the covariance with the \emph{bin-averaged} profile
$\bar P_{0,s}(y)$, the mean of $P_0(sp/N)$ over the (class, good prime) pairs of the bin of $y$ with the pair weights of the covariance;
the point value $P_0(sy)$ at the bin centre is given next to it (Table~\ref{tab:dil};
\texttt{anc/scripts/s70\_dilation\_fixedbins.py}).

\begin{figure}[ht]
\centering
\includegraphics[width=0.82\textwidth]{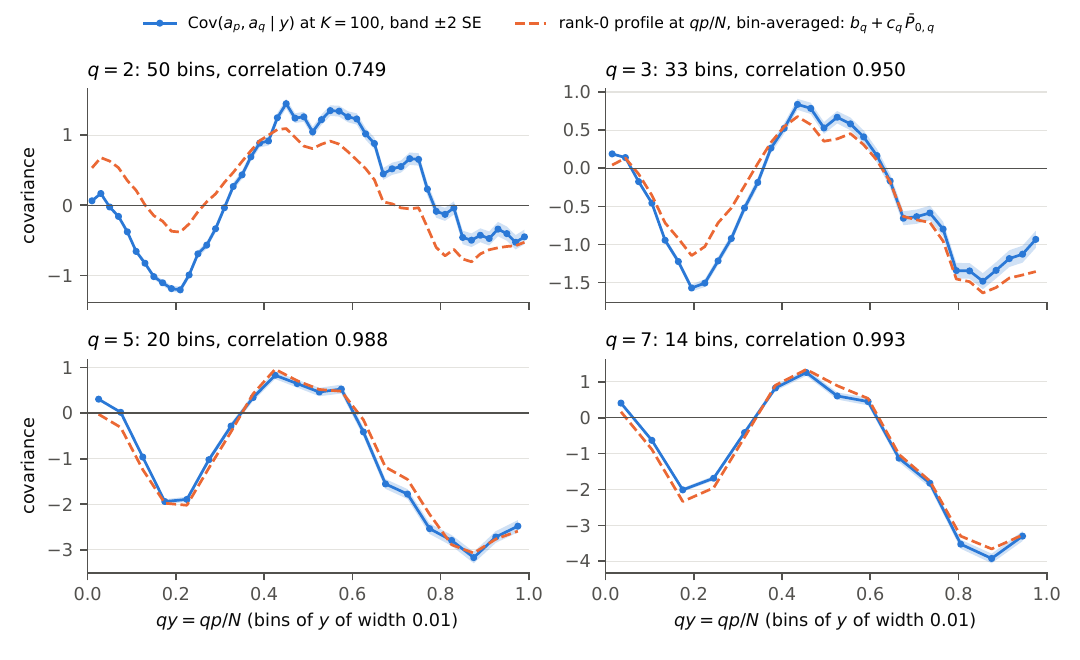}
\caption{$\Cov(\ap,a_q\mid y)$ for $q=2,3,5,7$ (solid; band $\pm2$ SE) as in Table~\ref{tab:dil}, against $qy$ over the bins with
$qy\le1$. Dashed: the bin-averaged profile $\bar P_{0,q}$ (not the point value $P_0(qy)$) as its affine least-squares fit
$b_q+c_q\bar P_{0,q}$ to the covariance; panel titles: the bin-averaged correlation of Table~\ref{tab:dil}(a). Source:
\texttt{anc/csv/dilation\_binavg.csv} (variant \texttt{all}; \texttt{s70\_dilation\_fixedbins.py}); figure \texttt{f02\_dilation.py}.}\label{fig:dilation}
\end{figure}

At the prespecified $s=q$, over the bins with $qy\le1$ (bin centre), $\Cov(\ap,a_q\mid y)$ and $\bar P_{0,q}(y)$ have correlation
$0.749$ ($q=2$, $50$ bins), $0.950$ ($q=3$, $33$ bins), $0.988$ ($q=5$, $20$ bins), $0.993$ ($q=7$, $14$ bins) and $0.999$ ($q=11$,
$9$ bins) (Figure~\ref{fig:dilation}); with the point value they are $0.749$, $0.941$, $0.988$, $0.994$ and $0.996$. Computing $P_0$
and the covariance on disjoint halves of the classes (classes paired consecutively in conductor order, one of each pair to each half at
random) changes the point-evaluation correlations at $s=q$ by at most $0.0045$ (\texttt{anc/scripts/s71\_split\_sample.py}). The
agreement is one of shape, and it is approximate: the residual of the affine fit $b_q+c_q\bar P_{0,q}$ exceeds the sampling noise, with
$\chi^2/\mathrm{dof}$ from $4.8$ ($q=11$) to $936$ ($q=2$) (Table~\ref{tab:dil}(a); per-bin variance $\mathrm{SE}_{\Cov}^2+c_q^2\,
\mathrm{SE}_{\bar P}^2$, with $\mathrm{SE}_{\bar P}$ (and its point analogue) the profile's per-bin SEs combined with the evaluation's weights, an upper bound, so that $\chi^2$ is if anything understated; correlations ignored), and from $14$ to $938$ with the point value (\texttt{s72\_residual\_noise.py}).

To compare dilation factors we scan $s$ over $[0.8q,1.2q]$ in steps of $0.05$ on one set of bins per $q$, the bins with $1.2q\,y\le1$
for every $y$ in the bin, so that all $s$ of a window are compared on the same part of the curves (Table~\ref{tab:dil}(b)). The correlation with $\bar P_{0,s}$ is largest at $s=1.75$, $2.90$, $5.10$, $7.10$ and $11.05$ for $q=2$, $3$,
$5$, $7$, $11$. For $q=5$, $7$, $11$ it is within $0.01$ of its maximum at $s=q$ and at least $0.97$ on $[5.00,5.25]$, $[6.80,7.35]$
and $[10.55,11.55]$; within these ranges the data do not single out $q$. For $q=3$ the maximum is $0.982$ at $s=2.90$, and $s=3$ gives
$0.964$, outside the set $[2.85,2.95]$ within $0.01$ of the maximum (point value: $0.955$, also outside). For $q=2$ the correlation is $0.743$ at $s=2$ against $0.866$ at $s=1.75$. The point evaluation gives the
same picture, and without the bin $y\le0.01$ no maximiser moves and each set changes by at most one grid step. \emph{Out of window:} in
the conductor windows $10^4\le N\le5\cdot10^4$ and $5\cdot10^4<N<10^5$ of Section~\ref{sec:profiles} (covariance, profile and pairs
from the classes of the window; $N<10^4$ is in neither) the maximisers agree to within $0.05$, and the $s$ chosen in one window gives in
the other a correlation within $0.004$ of that window's maximum (within $0.007$ on good $q$). The departures of the maximiser from $q$ at
$q=2$ and $3$ are thus the same in both windows; we do not interpret them.

Many classes have bad $q$: $q\mid N$ for $79.8\%$, $66.6\%$, $43.8\%$, $30.6\%$ and $17.7\%$ of the rank-$0$ class representatives for
$q=2$, $3$, $5$, $7$, $11$ (and for $80.1\%$ and $67.2\%$ of their pairs in the bins used, for $q=2$, $3$). For $q=2$ and $3$ the
statistic is therefore mostly over curves with bad $q$, where $a_q\in\{0,\pm1\}$ records the reduction type at $q$. On the classes with
$q\nmid N$ (variant \texttt{good\_q}: covariance, profile and pairs over these classes) the bin-averaged correlation at $s=q$ is
$0.508$ ($0.784$ against the profile of all rank-$0$ classes), $0.942$, $0.982$, $0.991$ and $0.998$, and for $q\ge3$ the maximisers
are within $0.05$ of those of all classes, but $s=q$ lies outside the set within $0.01$ of the maximum for $q=3$ and $q=5$ ($0.969$
against $0.982$ at $s=5.15$); for $q=2$ the correlation on the fixed bins increases towards the lower end of the window
($0.632$ at $s=1.6$, $0.380$ at $s=2$), which thus contains no interior maximum. The shortfall at $q=2$ thus persists on good $q$, and it is
unexplained.

In the affine fit $b_q+c_q\bar P_{0,q}$ (Figure~\ref{fig:dilation}), $b_q=+0.197$, $-0.425$, $-0.939$, $-0.999$, $-0.848$ for $q=2$,
$3$, $5$, $7$, $11$ and $c_q$ is in Table~\ref{tab:dil}(a); a correlation says nothing about the amplitude $c_q$, which we do not predict.
This regularity is empirical; we do not derive it.

\begin{table}[ht]
\centering
\caption{$\Cov(\ap,a_q\mid y)$ against the rank-$0$ profile dilated by $s$; rank-$0$ class representatives, $K=100$, good $p$ (pair
$p=q$ included), stored $a_q$ when $q\mid N$. $\rho$: correlation over the bins; ``bin-av.'': with $\bar P_{0,s}(y)$, the mean of
$P_0(sp/N)$ over the pairs of the bin; ``point'': with $P_0(sy)$ at the bin centre. (a) $s=q$, bins with $qy\le1$: $\rho$ (all classes,
bin-averaged and point; good $q$, bin-averaged), the amplitude $c_q$ of the affine least-squares fit and $\chi^2/\mathrm{dof}$ of its
residual. (b) $s=0.8q,0.8q+0.05,\dots,1.2q$ on the bins with $1.2q\,y\le1$ for every $y$ of the bin: maximiser $s_{\max}$ (first on the
grid; point evaluation in parentheses), $\rho_{\max}$, $\rho$ at $s=q$, the sets of $s$ with $\rho\ge\rho_{\max}-0.01$ and with
$\rho\ge0.97$; $s_{\max}$ on good $q$ and in the windows $W_1=[10^4,5\cdot10^4]$, $W_2=(5\cdot10^4,10^5)$ of $N$; ${}^{\dagger}$: end
of the window. Source: \texttt{anc/csv/dilation\_fixedbins.json} (\texttt{s70}); table \texttt{t06\_dilation\_table.py}.}\label{tab:dil}
{\small
\textup{(a)}\quad\begin{tabular}{r r ccc cc cc}
\toprule
 & & \multicolumn{3}{c}{correlation} & \multicolumn{2}{c}{$c_q$} & \multicolumn{2}{c}{$\chi^2/\mathrm{dof}$}\\
\cmidrule(lr){3-5}\cmidrule(lr){6-7}\cmidrule(lr){8-9}
$q$ & bins & bin-av. & point & good $q$ & bin-av. & point & bin-av. & point\\
\midrule
$2$ & $50$ & $0.749$ & $0.749$ & $0.508$ & $0.269$ & $0.269$ & $936.2$ & $937.5$\\
$3$ & $33$ & $0.950$ & $0.941$ & $0.942$ & $0.333$ & $0.322$ & $80.7$ & $88.6$\\
$5$ & $20$ & $0.988$ & $0.988$ & $0.982$ & $0.598$ & $0.586$ & $99.9$ & $46.6$\\
$7$ & $14$ & $0.993$ & $0.994$ & $0.991$ & $0.769$ & $0.718$ & $48.2$ & $17.4$\\
$11$ & $9$ & $0.999$ & $0.996$ & $0.998$ & $0.876$ & $0.861$ & $4.8$ & $13.7$\\
\bottomrule
\end{tabular}

\medskip
\textup{(b)}\quad\begin{tabular}{r r c cc c c c c}
\toprule
 & & \multicolumn{5}{c}{all rank-$0$ classes, bin-averaged} & good $q$ & $W_1$, $W_2$\\
\cmidrule(lr){3-7}
$q$ & bins & $s_{\max}$ (point) & $\rho_{\max}$ & $\rho(q)$ (point) & $\rho\ge\rho_{\max}-0.01$ & $\rho\ge0.97$ & $s_{\max}$ & $s_{\max}$\\
\midrule
$2$ & $41$ & $1.75$ (1.80) & $0.866$ & $0.743$ (0.742) & $[1.75,1.85]$ & --- & $1.60$$^{\dagger}$ & $1.80$, $1.75$\\
$3$ & $27$ & $2.90$ (2.95) & $0.982$ & $0.964$ (0.955) & $[2.85,2.95]$ & $[2.85,2.95]$ & $2.90$ & $2.90$, $2.90$\\
$5$ & $16$ & $5.10$ (5.05) & $0.985$ & $0.978$ (0.978) & $[5.00,5.25]$ & $[5.00,5.25]$ & $5.15$ & $5.10$, $5.10$\\
$7$ & $11$ & $7.10$ (7.00) & $0.993$ & $0.991$ (0.993) & $[6.95,7.25]$ & $[6.80,7.35]$ & $7.15$ & $7.05$, $7.10$\\
$11$ & $7$ & $11.05$ (11.05) & $0.998$ & $0.998$ (0.991) & $[10.80,11.30]$ & $[10.55,11.55]$ & $11.10$ & $11.10$, $11.05$\\
\bottomrule
\end{tabular}
}
\end{table}

(c) \emph{What Remark~\ref{rem:series} does not cover.} The quartile gaps of Section~\ref{sec:L} are nonlinear functions of
$L(E,1)$, and Remark~\ref{rem:series}, which is linear in $L(E,1)$, does not apply to them.

\subsection{Robustness}\label{sec:robust}
We checked four compositional or definitional alternatives to a dependence on $L(E,1)$; $K=40$ unless stated.
\begin{itemize}[leftmargin=1.5em]
\item \emph{Conductor.} Forming the quartiles within $20$ conductor bins instead of over all $N$ changes $Q_4-Q_1$ by at most $0.061$ at
$y\le0.125$ (computed on the $K=100$ bin sums; \texttt{pn\_check\_condbins.py}).
\item \emph{Reduction type.} Standardising each quartile's profile to the rank-$0$ distribution of the pair $(m,a)$, where
$m=\min(\#\text{multiplicative primes},5)$ and $a=\min(\#\text{additive primes},3)$ (quartiles over all $N$), changes the gap by at most
$0.31$ over the $40$ bins; the largest shifts are on $0.9<y\le1$ ($0.25$ to $0.31$; at $y\approx0.96$ from $+1.49$ to $+1.80$), and at
$y\approx0.0625$ the gap moves from $-3.822$ to $-3.806$. No bin changes sign, in particular none of the $31$ bins with $|z|\ge3.23$
(\texttt{anc/scripts/s61\_robust\_allbins.py}).
\item \emph{Drift toward rank $2$.} With $P_r$ the rank-$r$ profile, a mixture $P_0+\lambda(P_2-P_0)$ (quartiles over all $N$; $40$
bins) does not describe the quartile profiles: with the standard error of each bin's residual $Q_j-(1-\lambda)P_0-\lambda P_2$
(linearised and clustered by class, so that it includes the sampling errors of $P_0$ and $P_2$ and the overlap of $Q_j$ with the
rank-$0$ sample) and $\lambda$ chosen to minimise this $\chi^2$, whose standard errors depend on $\lambda$, $\chi^2/\mathrm{dof}$ lies between $231$ and $980$ (\texttt{anc/scripts/s45\_rank2\_chi2.py}). This
rules out only that the quartile profiles are mixtures of this form.
\item \emph{Bad primes.} The profiles use good primes. Including the stored entries at $p\mid N$ moves the lowest $K=100$ bin of the
rank-$0$ profile by $10$ cluster SE ($1.018$ to $0.979$) and the next five bins ($0.01<y\le0.06$) by less than $0.2$ SE.
\end{itemize}
The dependence on $L(E,1)$ survives the conductor and reduction-type adjustments tested.

\section{Weighting by BSD invariants}\label{sec:weights}

Over all ranks, with $\varepsilon_E=(-1)^{\mathrm{rank}}$, let $D(y)=\sum_E\varepsilon_E S_{1,E}/\sum_E n_E$ in each bin (class
representatives, good primes, $K=40$), and for a weight $w$ let $D_w(y)=\sum_E w_E\varepsilon_E S_{1,E}/\sum_E w_E n_E$, with the
representative's value of the invariant as $w_E$. Table~\ref{tab:weights} gives $D_w-D$ at $y\le0.15$; $z$ uses a bootstrap over classes
($B=1000$).

\begin{itemize}[leftmargin=1.5em]
\item The $\Tam$ weight changes $D$ by $+0.03$ ($z=+2.8$) at $y\le0.025$ and by $-0.89$, $-1.21$, $-1.24$, $-1.39$ ($z=-31$, $-33$,
$-28$, $-22$) on the next four bins.
\item With class-level weights (the mean, maximum or sum of $\Tam$ over the curves of the class) the drop on these four bins is between
$4\%$ and $21\%$ smaller and has the same sign.
\item On these four bins the change after fixing the reduction-type composition $(m,a)$ of Section~\ref{sec:robust} is about $70\%$ of
the change before, after fixing the rank composition $89$--$96\%$, and after fixing the conductor quantiles $99\%$ (direct
standardisation to the unweighted strata shares). At $y\le0.025$ the change is small and reverses sign under the type fixing ($+0.035$,
$z=+2.8$, becomes $-0.078$, $z=-13$; Table~\ref{tab:weights}), so no such ratio is defined there.
\item The $\Omeg$ weight (all ranks) changes $D$ by $-0.18$, $-0.18$, $-0.19$, $-0.13$ at $y\le0.1$.
\end{itemize}

On rank $0$ the $L(E,1)$ weight changes the mean of $\ap$ by $-0.22$, $-0.81$, $-0.96$, $-0.95$ at $y\le0.1$ ($z$ down to
$-117$), and the $\Tam$ weight by $-0.11$, $-0.75$, $-1.21$, $-1.24$ (Table~\ref{tab:weights}). The two invariants are correlated,
$\Corr(\log\Tam,\log L(E,1))=0.33$ on the rank-$0$ class representatives, and we compare the weights by direct standardisation: the
weighted mean is computed within strata and recombined with the unweighted pair shares of the strata (Tables~\ref{tab:std} and~\ref{tab:stdb}; standard errors
from a bootstrap over classes, $B=1000$; \texttt{anc/scripts/s31\_weights\_v4.py}). For any weight $w(E)>0$ the change of the rank-$0$
mean is exactly $M_b(w\ap)/M_b(w)-M_b(\ap)=\Cov_b(\ap,w)/M_b(w)$, so the $L(E,1)$-weighted change is $\Cov_b(\ap,L(E,1))/M_b(L(E,1))$,
which \eqref{eq:split} rewrites. Below, bin $j$ is $(j/40,(j+1)/40]$, ``std.'' means standardised to the deciles of $\log L(E,1)$, and
$\pm$ is one bootstrap SE. The percentages below are ratios of the standardised to the unstandardised profile change (or of their
sums over bins), not shares of a separately identified arithmetic effect of the weight.
\begin{itemize}[leftmargin=1.5em]
\item $\Tam$: the weight is concentrated (effective number of classes $5{,}869$ of $168{,}760$; the top $1\%$ of classes carry $31\%$ of it;
largest value $13{,}824$). On the bins $1$--$4$ the std.\ change (Table~\ref{tab:std}) is $-11\%$, $9\%$, $15\%$ and $30\%$ of the
unstandardised change; summed over these bins it is $13.2\%\pm1.4\%$ of the summed unstandardised change,
and the sum of its absolute values is $16.9\%\pm2.2\%$ of theirs. In bin $0$ it is $+0.125\pm0.013$ against $-0.106\pm0.024$ unstandardised,
opposite in sign and larger in absolute value; in bin $5$ it is $-0.047\pm0.069$ against $-0.820\pm0.119$.
\item $\Omeg$: on the bins $1$--$5$ the unstandardised change lies between $-0.259$ and $-0.183$ (SE $0.008$ to $0.017$; Table~\ref{tab:std});
on the bins $3$--$5$ the std.\ change has the opposite sign ($12\%$, $28\%$, $25\%$ of the unstandardised size). Over the bins $1$--$4$ the signed
ratio is $-6.1\%\pm2.1\%$ and the ratio of absolute values $13.7\%\pm2.4\%$. In bin $0$: $+0.007\pm0.003$ against $+0.037\pm0.004$.
\item $L(E,1)$, standardised to the Tamagawa classes $\Tam=1$, $2$, $3$--$4$, $5$--$8$, $9$--$16$, $\ge17$: the standardised change is $83\%$, $91\%$, $89\%$,
$87\%$, $84\%$, $87\%$ of the unstandardised one on the bins $0$--$5$ (SE $0.3$ to $1.0$ percentage points).
\end{itemize}
With $10$, $20$, $50$, $100$, $200$ and $500$ quantile strata of $\log L(E,1)$ (Table~\ref{tab:stdb}; \texttt{s43\_strata\_stress.py},
$B=1000$) the std.\ $\Tam$ results persist: the signed ratio over the bins $1$--$4$ lies
between $11.5\%$ and $13.3\%$, the ratio of absolute values between $15.0\%$ and $16.9\%$, the bin-$0$ value between $+0.123$ and $+0.131$, and no
bin moves by more than $0.066$ from its decile value (bin $4$, from $-0.380$ to $-0.314$, with $z$ from $-5.8$ to $-7.9$). The std.\
$\Omeg$ change on the bins $4$ and $5$ is positive at every stratification tried. It decreases as the strata are refined, from $+0.058$
and $+0.055$ with deciles to $+0.040$ and $+0.033$ with $100$ strata ($z=3.0$ and $2.1$), and levels off there ($+0.039$ and $+0.031$
with $200$ and with $500$ strata, $z=2.8$--$2.9$ and $2.0$); in bin $3$ it falls from $+0.029\pm0.012$ to $+0.008\pm0.011$
($500$ strata). No adjustment was made for the
multiplicity of bins, weights and stratifications compared. The signed $\Omeg$ ratio over the bins $1$--$4$ changes sign with the
stratification: $-6.1\%$, $-2.0\%$, $+1.5\%$, $+2.8\%$, $+3.6\%$ and $+3.7\%$ for $10$ to $500$ strata (ratio of absolute values $13.7\%$ to
$14.7\%$). The $\Tors^2$ weight changes the rank-$0$ mean by at most $0.113$ in absolute value on the bins $0$--$5$, and by at most $0.138$ after
standardisation (Table~\ref{tab:std}). These standardisations are on
rank $0$ only: at odd rank $L(E,1)=0$, and the all-ranks root-number-weighted $D$ above is not standardised in this way; for it the source of
the changes is open. The same analysis at $p\le3571$ on all curves is in \texttt{anc/p3571\_weights/}.

\begin{table}[ht]
\centering
\caption{Weighted root-number sums at $K=40$, bins $y\in(j/40,(j+1)/40]$, $j=0,\dots,5$; class representatives ($437{,}226$; rank $0$:
$168{,}760$), good primes. $D(y)$: all ranks, $\varepsilon=(-1)^{\mathrm{rank}}$; ``mean $\ap$'': rank $0$. $\Delta_w=D_w-D$ for the weight
$w$ (representative's value); below it, in parentheses, $z=\Delta_w/\mathrm{SE}$, SE from a bootstrap over classes ($B=1000$, seed $20261001$).
``Type-fixed'': the $\Tam$ change after direct standardisation of the within-stratum weighted $D$ to the unweighted pair shares of the
strata $(m,a)$ of Section~\ref{sec:robust}, bootstrap SE ($B=1000$, same seed). $|z|\ge10$ is rounded to an integer. Sources:
\texttt{data\_pN/results/pnB\_weighted\_K40.csv}, \texttt{anc/csv/typefixed\_boot.csv} (script \texttt{s48\_typefixed\_boot.py}); table
\texttt{t01\_tables.py}.}\label{tab:weights}
{\small\renewcommand{\arraystretch}{1.0}
\begin{tabular}{l r rr r @{\hspace{1.6em}} r rr}
\toprule
 & \multicolumn{4}{c}{all ranks, $\varepsilon=(-1)^{\mathrm{rank}}$} & \multicolumn{3}{c}{rank $0$}\\
\cmidrule(lr){2-5}\cmidrule(lr){6-8}
$y$ & $D(y)$ & $\Delta_{\prod c_p}$ & $\Delta_{\prod c_p}$ & $\Delta_{\Omega}$ & mean $a_p$ & $\Delta_{L(E,1)}$ & $\Delta_{\prod c_p}$\\
 & & & type-fixed & & & & \\
\midrule
$(0,1/40]$ & $+0.832$ & $+0.035$ & $-0.078$ & $-0.178$ & $+1.310$ & $-0.222$ & $-0.106$\\
 & & $(+2.8)$ & $(-13)$ & $(-71)$ & & $(-69)$ & $(-4.7)$\\[1pt]
$(1/40,2/40]$ & $+1.462$ & $-0.891$ & $-0.636$ & $-0.176$ & $+1.798$ & $-0.812$ & $-0.750$\\
 & & $(-31)$ & $(-45)$ & $(-38)$ & & $(-117)$ & $(-14)$\\[1pt]
$(2/40,3/40]$ & $+1.151$ & $-1.206$ & $-0.889$ & $-0.190$ & $+1.434$ & $-0.958$ & $-1.214$\\
 & & $(-33)$ & $(-50)$ & $(-30)$ & & $(-104)$ & $(-13)$\\[1pt]
$(3/40,4/40]$ & $+0.434$ & $-1.241$ & $-0.891$ & $-0.126$ & $+0.637$ & $-0.952$ & $-1.244$\\
 & & $(-28)$ & $(-42)$ & $(-17)$ & & $(-85)$ & $(-12)$\\[1pt]
$(4/40,5/40]$ & $-0.345$ & $-1.394$ & $-0.951$ & $-0.161$ & $-0.035$ & $-0.745$ & $-1.283$\\
 & & $(-22)$ & $(-34)$ & $(-19)$ & & $(-57)$ & $(-13)$\\[1pt]
$(5/40,6/40]$ & $-1.116$ & $-0.319$ & $-0.265$ & $-0.123$ & $-0.979$ & $-0.543$ & $-0.820$\\
 & & $(-5.1)$ & $(-9.2)$ & $(-13)$ & & $(-37)$ & $(-7.1)$\\
\bottomrule
\end{tabular}
}
\end{table}

\begin{table}[ht]
\centering
\caption{Rank-$0$ class representatives ($168{,}760$), good primes, $K=40$, bins $y\in(j/40,(j+1)/40]$: change $\Delta_w$ of the
pairs-weighted mean of $\ap$ under the weight $w$, and the change after direct standardisation (within-stratum weighted means
recombined with the unweighted pair shares of the strata): ``std.\ ($L$)'' to the deciles of $\log L(E,1)$, ``std.\ (Tam.)'' to the
Tamagawa classes $\Tam=1$, $2$, $3$--$4$, $5$--$8$, $9$--$16$, $\ge17$. Below each value, in parentheses: bootstrap SE over classes
($B=1000$, seed $20261002$; each class keeps its stratum). Source: \texttt{anc/csv/weights\_v4.csv} (script \texttt{s31\_weights\_v4.py});
table \texttt{t05\_std\_table.py}.}\label{tab:std}
{\small\renewcommand{\arraystretch}{1.0}
\begin{tabular}{l rr rr rr rr}
\toprule
 & \multicolumn{2}{c}{$w=\prod c_p$} & \multicolumn{2}{c}{$w=\Omega_E$} & \multicolumn{2}{c}{$w=L(E,1)$} & \multicolumn{2}{c}{$w=|E(\mathbb Q)_{\mathrm{tors}}|^2$}\\
\cmidrule(lr){2-3}\cmidrule(lr){4-5}\cmidrule(lr){6-7}\cmidrule(lr){8-9}
$y$ & $\Delta_w$ & std.\ ($L$) & $\Delta_w$ & std.\ ($L$) & $\Delta_w$ & std.\ (Tam.) & $\Delta_w$ & std.\ ($L$)\\
\midrule
$(0/40,1/40]$ & $-0.106$ & $+0.125$ & $+0.037$ & $+0.007$ & $-0.222$ & $-0.185$ & $+0.030$ & $+0.029$\\
 & $(0.024)$ & $(0.013)$ & $(0.004)$ & $(0.003)$ & $(0.003)$ & $(0.003)$ & $(0.005)$ & $(0.005)$\\[1pt]
$(1/40,2/40]$ & $-0.750$ & $+0.082$ & $-0.183$ & $-0.028$ & $-0.812$ & $-0.736$ & $+0.069$ & $+0.050$\\
 & $(0.053)$ & $(0.035)$ & $(0.008)$ & $(0.007)$ & $(0.007)$ & $(0.007)$ & $(0.012)$ & $(0.011)$\\[1pt]
$(2/40,3/40]$ & $-1.214$ & $-0.109$ & $-0.259$ & $-0.005$ & $-0.958$ & $-0.849$ & $+0.032$ & $+0.003$\\
 & $(0.095)$ & $(0.052)$ & $(0.010)$ & $(0.009)$ & $(0.009)$ & $(0.009)$ & $(0.016)$ & $(0.014)$\\[1pt]
$(3/40,4/40]$ & $-1.244$ & $-0.187$ & $-0.238$ & $+0.029$ & $-0.952$ & $-0.825$ & $-0.015$ & $-0.044$\\
 & $(0.096)$ & $(0.070)$ & $(0.013)$ & $(0.012)$ & $(0.012)$ & $(0.011)$ & $(0.018)$ & $(0.017)$\\[1pt]
$(4/40,5/40]$ & $-1.283$ & $-0.380$ & $-0.203$ & $+0.058$ & $-0.745$ & $-0.629$ & $-0.113$ & $-0.138$\\
 & $(0.099)$ & $(0.065)$ & $(0.015)$ & $(0.014)$ & $(0.013)$ & $(0.012)$ & $(0.021)$ & $(0.020)$\\[1pt]
$(5/40,6/40]$ & $-0.820$ & $-0.047$ & $-0.216$ & $+0.055$ & $-0.543$ & $-0.470$ & $-0.112$ & $-0.134$\\
 & $(0.119)$ & $(0.069)$ & $(0.017)$ & $(0.016)$ & $(0.014)$ & $(0.014)$ & $(0.022)$ & $(0.022)$\\[1pt]
\bottomrule
\end{tabular}
}
\end{table}

\begin{table}[ht]
\centering
\caption{The standardised $\Tam$ and $\Omeg$ changes of Table~\ref{tab:std} with $10$ (deciles), $20$, $50$, $100$, $200$ and
$500$ quantile strata of $\log L(E,1)$; over the bins $1$--$4$, the ratio in percent of the standardised to the unstandardised change,
of their sums (signed) and of their sums of absolute values (absolute). SEs as in Table~\ref{tab:std}. Source:
\texttt{anc/csv/strata\_stress.csv}, \texttt{strata\_stress.json} (script \texttt{s43\_strata\_stress.py}); table
\texttt{t05\_std\_table.py}.}\label{tab:stdb}
{\small\renewcommand{\arraystretch}{1.0}
\begin{tabular}{l r rrrrrr rr}
\toprule
 & & \multicolumn{6}{c}{std.\ change in the bin $(j/40,(j+1)/40]$} & \multicolumn{2}{c}{bins $1$--$4$: ratio (\%)}\\
\cmidrule(lr){3-8}\cmidrule(lr){9-10}
$w$ & strata & $j=0$ & $j=1$ & $j=2$ & $j=3$ & $j=4$ & $j=5$ & signed & absolute\\
\midrule
$\prod c_p$ & $10$ & $+0.125$ & $+0.082$ & $-0.109$ & $-0.187$ & $-0.380$ & $-0.047$ & $+13.2$ & $16.9$\\
 &  & $(0.013)$ & $(0.035)$ & $(0.052)$ & $(0.070)$ & $(0.065)$ & $(0.069)$ & $(1.4)$ & $(2.2)$\\[1pt]
 & $20$ & $+0.131$ & $+0.082$ & $-0.110$ & $-0.186$ & $-0.383$ & $-0.061$ & $+13.3$ & $16.9$\\
 &  & $(0.013)$ & $(0.032)$ & $(0.048)$ & $(0.065)$ & $(0.063)$ & $(0.066)$ & $(1.3)$ & $(2.1)$\\[1pt]
 & $50$ & $+0.128$ & $+0.081$ & $-0.110$ & $-0.177$ & $-0.373$ & $-0.067$ & $+12.9$ & $16.5$\\
 &  & $(0.012)$ & $(0.029)$ & $(0.043)$ & $(0.056)$ & $(0.058)$ & $(0.061)$ & $(1.2)$ & $(1.8)$\\[1pt]
 & $100$ & $+0.127$ & $+0.080$ & $-0.105$ & $-0.173$ & $-0.354$ & $-0.068$ & $+12.3$ & $15.9$\\
 &  & $(0.011)$ & $(0.026)$ & $(0.038)$ & $(0.050)$ & $(0.052)$ & $(0.056)$ & $(1.1)$ & $(1.7)$\\[1pt]
 & $200$ & $+0.127$ & $+0.079$ & $-0.103$ & $-0.181$ & $-0.335$ & $-0.080$ & $+12.0$ & $15.5$\\
 &  & $(0.010)$ & $(0.023)$ & $(0.033)$ & $(0.044)$ & $(0.047)$ & $(0.051)$ & $(1.1)$ & $(1.5)$\\[1pt]
 & $500$ & $+0.123$ & $+0.078$ & $-0.098$ & $-0.183$ & $-0.314$ & $-0.092$ & $+11.5$ & $15.0$\\
 &  & $(0.008)$ & $(0.019)$ & $(0.028)$ & $(0.036)$ & $(0.040)$ & $(0.043)$ & $(1.0)$ & $(1.3)$\\[1pt]
\midrule
$\Omega_E$ & $10$ & $+0.007$ & $-0.028$ & $-0.005$ & $+0.029$ & $+0.058$ & $+0.055$ & $-6.1$ & $13.7$\\
 &  & $(0.003)$ & $(0.007)$ & $(0.009)$ & $(0.012)$ & $(0.014)$ & $(0.016)$ & $(2.1)$ & $(2.4)$\\[1pt]
 & $20$ & $-0.006$ & $-0.041$ & $-0.014$ & $+0.020$ & $+0.052$ & $+0.046$ & $-2.0$ & $14.4$\\
 &  & $(0.003)$ & $(0.007)$ & $(0.009)$ & $(0.012)$ & $(0.014)$ & $(0.016)$ & $(2.0)$ & $(2.5)$\\[1pt]
 & $50$ & $-0.012$ & $-0.048$ & $-0.024$ & $+0.015$ & $+0.043$ & $+0.036$ & $+1.5$ & $14.7$\\
 &  & $(0.003)$ & $(0.007)$ & $(0.009)$ & $(0.012)$ & $(0.014)$ & $(0.016)$ & $(1.9)$ & $(2.5)$\\[1pt]
 & $100$ & $-0.014$ & $-0.051$ & $-0.026$ & $+0.011$ & $+0.040$ & $+0.033$ & $+2.8$ & $14.6$\\
 &  & $(0.003)$ & $(0.007)$ & $(0.009)$ & $(0.012)$ & $(0.014)$ & $(0.016)$ & $(1.9)$ & $(2.4)$\\[1pt]
 & $200$ & $-0.014$ & $-0.053$ & $-0.028$ & $+0.010$ & $+0.039$ & $+0.031$ & $+3.6$ & $14.7$\\
 &  & $(0.003)$ & $(0.007)$ & $(0.009)$ & $(0.011)$ & $(0.014)$ & $(0.016)$ & $(1.9)$ & $(2.4)$\\[1pt]
 & $500$ & $-0.015$ & $-0.052$ & $-0.029$ & $+0.008$ & $+0.039$ & $+0.031$ & $+3.7$ & $14.6$\\
 &  & $(0.003)$ & $(0.007)$ & $(0.009)$ & $(0.011)$ & $(0.014)$ & $(0.015)$ & $(1.9)$ & $(2.4)$\\[1pt]
\bottomrule
\end{tabular}
}
\end{table}

\section{Prime conductors and the Sawin--Sutherland prediction}\label{sec:ss}

Sawin and Sutherland \cite[\S1.5]{SS} state, for curves of prime conductor ordered by conductor, a formula for the limit as $X\to\infty$
of the root-number-weighted sums of $\ap$ over $p\in(C_1N,C_2N]$, their equation~(7) (their \S1.5, (7), (8) and Figure~4 are in v2 of the
preprint), introduced by ``It seems likely that for real
numbers $0<C_1<C_2$ we have''. We use (7) and (8) exactly as displayed in v2, with the local factor $\tilde\ell'_{p,2v_p(m)}$ only at
the primes $p\mid m$. Their sentence after (8), that (7) agrees with their Conjecture~1 except that it uses ``the local factors
$\tilde\ell'_{p,\nu}$ instead of $\ell_{p,\nu}$ and $\hat\ell_{p,\nu}$'', admits a reading with further factors at primes $p\mid q$,
$p\nmid m$; the rejection statements below are conditional on the displayed reading. Of the status of (7) they write: ``Since the reasoning is non-rigorous in multiple steps, we do not give any proofs
and simply state a prediction, an informal justification, and give numerical evidence for it''. Their Figure~4 compares (7) with curves
of prime conductor up to $X=2^{30}$.

\begin{figure}[ht]
\centering
\includegraphics[width=0.88\textwidth]{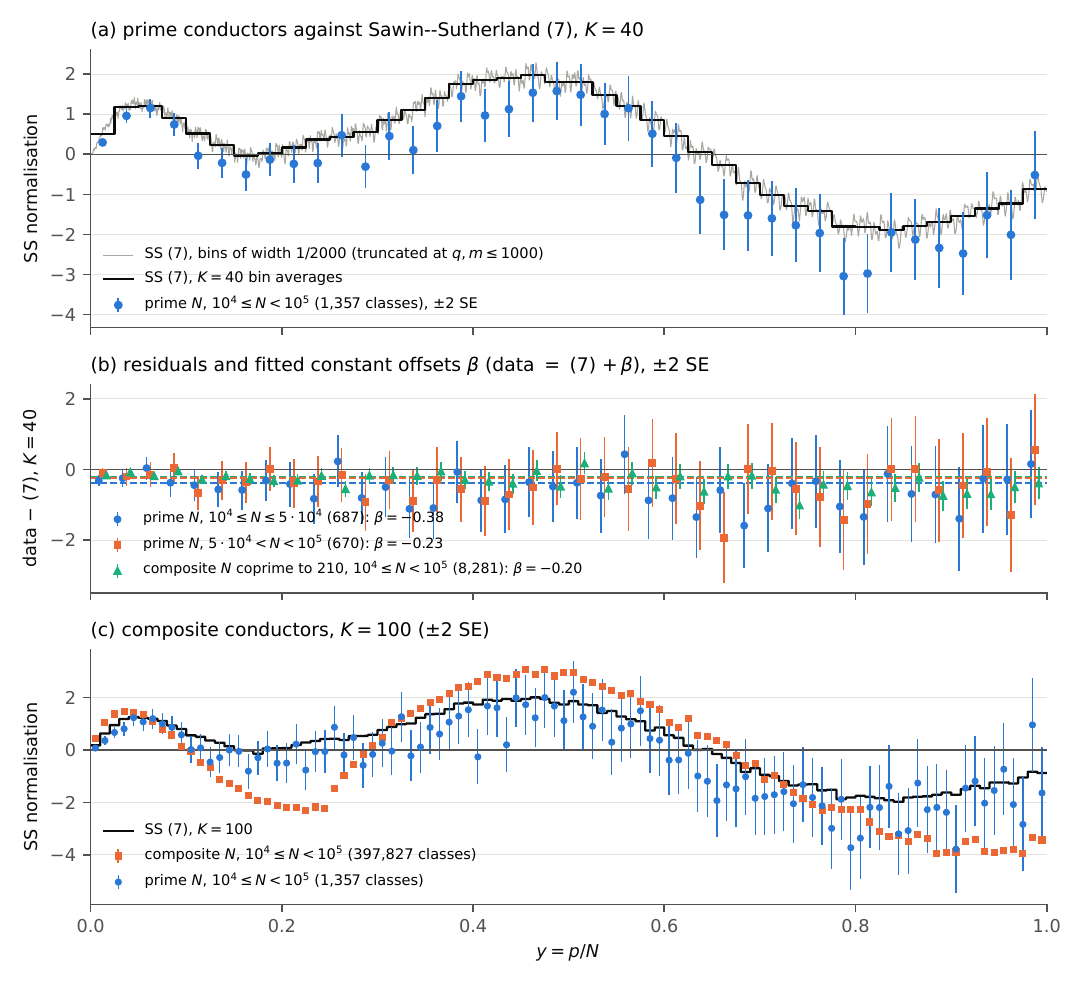}
\caption{Sawin--Sutherland (7) \cite[\S1.5]{SS} and prime-conductor classes, $10^4\le N<10^5$. Data $D_j$ in the normalisation of (7):
$\varepsilon=(-1)^{\mathrm{rank}}$, equal weight per class, all primes of the bin \emph{including} $p\mid N$ (unlike the other figures,
which use good primes only); bars $\pm2$ curve-cluster SE. Prediction: bin averages of the right-hand side of (7), sums truncated at
squarefree $q\le3000$, $m\le10^4$; grey, bins of width $1/2000$ truncated at $q,m\le1000$. (a) $K=40$. (b) Data minus prediction in two
conductor windows and for composite $N$ coprime to $210$, with the fitted offsets $\beta$ (dashed). (c) $K=100$: composite $N$
($397{,}827$ classes) and prime $N$. Sources: \texttt{anc/csv/ss\_compare\_K40.csv}, \texttt{ss\_compare\_K100.csv}, \texttt{ss\_summary.csv},
\texttt{ss\_prediction\_fine.csv}; script \texttt{anc/scripts/f01\_figures.py}.}\label{fig:ss}
\end{figure}

We compare (7) with the $1{,}357$ prime-conductor classes with $10^4\le N<10^5$ (Figure~\ref{fig:ss}), in the
normalisation of (7): in bin $j$, $D_j$ is $K$ times the mean over classes of $(\log(N(j+\tfrac12)/K)/N)\sum_{p}\varepsilon\,\ap$, the sum
over all primes of the bin including $p\mid N$; the bin averages $P_j$ of the right-hand side of (7) are computed with the sums truncated
at squarefree $q\le3000$ and $m\le10^4$. Equation (7) predicts a limit as $X\to\infty$, obtained by an approach that, in Sawin and
Sutherland's words, ``relies on multiple heuristics'' \cite[\S1.5]{SS}; a discrepancy at $N<10^5$ cannot test that limit, and each
rejection below is a statement about this finite range only. At $K=40$ the data for $10^4\le N<10^5$ do not fit (7) as stated
($\chi^2/\mathrm{dof}=3.64$, diagonal errors). A pure amplitude fit $D=\alpha P$ gives
$\alpha=0.82\pm0.04$ and does not fit these data either ($\chi^2/\mathrm{dof}=120.3/39=3.09$; \texttt{anc/csv/ss\_summary.csv}). A constant offset, $D=P+\beta$, gives $\beta=-0.31\pm0.03$ with
$\chi^2/\mathrm{dof}=1.43$ ($p=0.04$) with diagonal errors. The bins are correlated: with the full covariance matrix of the per-class
vectors the offset fit gives $\beta=-0.32\pm0.03$ and $\chi^2=67.5$ on $39$ degrees of freedom ($p=0.003$; $p=0.0055$ when $\chi^2$ is
calibrated by a bootstrap over classes, $B=2000$; \texttt{anc/csv/ss\_fullcov.json}), so on the pooled window the offset model is
rejected as well, for these finite-$N$ data (not a test of the limit). The $40\times40$ covariance is estimated from $670$ to $1{,}686$ classes, so the nominal $\chi^2$ $p$-values of the
full-covariance fits are mildly anti-conservative; calibrated by Hotelling's $T^2$ distribution ($F$ with $39$ and $n-39$ degrees of
freedom for $n$ classes) they are $0.0057$ (pooled), $0.11$ and $0.46$ (the two windows below) and $0.0048$ (all prime $N<10^5$), and no
conclusion of this section changes (\texttt{anc/scripts/s62\_ss\_hotelling.py}). Within each of the two conductor windows the offset model is not rejected at the $5\%$ level: $\beta=-0.38\pm0.05$ ($10^4\le N\le5\cdot10^4$, $687$
classes) and $-0.23\pm0.05$ ($5\cdot10^4<N<10^5$, $670$ classes), a difference of $2.2\sigma$, with $p=0.33$ and $0.59$ for the offset fit
with diagonal errors and $p=0.060$ and $0.35$ with the full covariance (\texttt{anc/scripts/s44\_ss\_windows.py}). The residual shapes of
the two windows are positively correlated ($0.26$ over the $40$ bins) and their difference is consistent with noise ($\chi^2=22.6$ on
$39$ degrees of freedom), so the pooled rejection reflects the larger sample: a misfit too small to reach the $5\%$ level in either window alone (first window,
full covariance: $p=0.060$). (7) has no
constant term. For curves ordered by height, Sawin and Sutherland observe ``a persistent downward bias, where (1) is slightly less than its
predicted limiting value'', which ``seems to decay as $X$ goes to $\infty$''. They raise curves of rank $\ge2$ as ``a possible
explanation'' and then argue against it: over their range of $X$ the proportion of rank-$2$ curves increases while the bias decreases,
the opposite of what a rank-$2$ source would lead one to expect \cite[\S1]{SS}. Of our
prime-conductor classes, $23\%$ have rank $\ge2$, against $10.5\%$ of all classes with $N<10^5$. (7) is a statement about the average
over all prime-conductor curves and makes no prediction for a subset selected by rank. The rank-$0$/$1$ fit is descriptive: on the
$1{,}040$ classes of rank $0$ and $1$ alone the offset is $\beta=+0.245\pm0.033$ with diagonal errors ($\chi^2=102.0$ on $39$ degrees of
freedom, $p=2\cdot10^{-7}$) and $\beta=+0.213\pm0.021$ with the full covariance ($\chi^2=114.4$, $p=2\cdot10^{-9}$), of the opposite sign
to the all-ranks value $-0.31\pm0.03$ ($-0.32\pm0.03$) (\texttt{anc/csv/ss\_rank01.json}). The average in (7) is over curves, ours over isogeny classes ($1{,}357$ classes, $1{,}383$ curves); weighting by curves gives
$\beta=-0.30\pm0.03$. Over all prime $N<10^5$ ($1{,}686$ classes) the offset fit has $p=0.08$ with diagonal errors and $p=0.003$
with the full covariance; this too is a small-$N$ check of a limit prediction and cannot test the limit. Composite conductors follow a
different curve (Figure~\ref{fig:ss}(c)); (7) is stated for prime conductors only. The difference
between composite $N$ coprime to $210$ and prime $N$ is not significant at this size ($\chi^2/\mathrm{dof}=1.38$), but this reflects
limited power: these $8{,}281$ composite classes themselves differ from the right-hand side of (7) ($\chi^2=352.2$ on $40$ degrees of
freedom) and from it plus an offset ($\chi^2=118.6$ on $39$), with diagonal errors.

\section{Relation to prior work, open questions, and corrections to the earlier versions}\label{sec:prior}

\subsection{Prior work}\label{sec:lit}
The results nearest to Sections~\ref{sec:small} and~\ref{sec:scale} that we found are the following.
\begin{itemize}[leftmargin=1.5em]
\item \emph{Central values and a fixed prime.} Knightly and Reno \cite[Theorem~1.5, Remark~3.4]{KR} and, for the product
$L(\frac12,\varphi\otimes\chi)L(\frac12,\varphi)$ at prime level, Ramakrishnan and Rogawski \cite[Theorem~A]{RR}, \cite{RRe} (on its
interval version see \cite[Remark~5.3(1)]{KR}), both for weight $k>2$ (Section~\ref{sec:small}). Feigon
and Whitehouse obtain exact formulas for averages of central values of twisted base change $L$-functions over Hilbert modular forms and
``also establish an equidistribution result for the Hecke eigenvalues weighted by these $L$-values'' \cite[Abstract]{FW}. Fazzari studies the low-lying zeros of families of $L$-functions
``tilted by a power of the $L$-function evaluated at the central point'' \cite[Abstract]{Faz}.
\item \emph{Central values and $p\asymp N$.} Martin \cite[\S7]{Mar} considers weight-$2$ newforms $f$ of prime level $N$ inert in
$\Q(\sqrt{-3})$ and the averages of $a_p(f)$ weighted by $L(1/2,f)L(1/2,f\otimes\chi_{-3})/(f,f)$, summed over $f$ and averaged over the
levels $X\le N\le2X$, for $p\le X$. His (7.1) is not a formula for these level averages: it is an exact identity at each such level
$N$ for the sum over $f$ of $L(1/2,f)L(1/2,f\otimes\chi_{-D})a_n(f)/(f,f)$ (stated for $N\neq p$; as printed, one of its terms contains
$p$ in place of $n$), and his averages are formed from it. Of such exact formulas he writes ``This was carried out in [FW09, Mar22] under
some hypotheses on $N$ and $D$'', and the setting of (7.1) ``was originally considered in [MR12] (using the Gross--Zagier formula, rather
than Jacquet's relative trace formula)'' \cite[\S7]{Mar}, where [FW09] is \cite{FW} and [MR12] is \cite{MR}. There the identity is an instance of
Michel and Ramakrishnan's Theorem~1, which they state ``for any positive integer $m$'' as an ``exact identity'' \cite[Theorem~1]{MR},
and whose last sum vanishes in what they call ``the stable range'' $N>mD$
\cite[Corollary~1]{MR}. Of his graphs Martin writes: ``These
graphs are apparently scale invariant in $p/X$, and we expect the blue and red graphs to tend to a limiting function after appropriate
smoothing''; in the stable range $p<N/D$ the sum on the right of (7.1) vanishes ``(in this range, it is not hard to prove the limit
exists)''. His weight is a product of two central values with the harmonic weight $1/(f,f)$, his family consists of all newforms, and he
weights rather than conditions.
\item Bettin \cite[Theorem~1.1]{Bet}, for the harmonic average $M_m(\alpha,\chi;N)$ of $a_f(m)L(\tfrac12+\alpha,f\otimes\chi)$ over
the primitive forms of weight $2$ and level $N$ (his (2.1) writes the sum over the Hecke basis of $S_2(N)$, which for prime $N$ is the set of
primitive forms; for $N=p^\nu$, $\nu\ge2$, his \S5 reduces the primitive-form average to such sums at levels $N$ and $N/p$): ``Let $N=p^\nu$ with $p$ prime and $\nu\ge1$ and let $\chi$ be a
primitive character modulo $q$ with $(q,N)=1$. Let $|\Re(\alpha)|\ll\frac1{\log N}$ and write $T=1+|\Im(\alpha)|$. Then, if $\nu\ge2$ and
$p|m$ then $M_m(\alpha,\chi;N)=0$. In all other cases for all $\varepsilon>0$ we have $M_m(\alpha,\chi;N)=\frac{\chi(m)}{m^{\frac12+\alpha}}
(1-\delta_\nu(p))+O_\varepsilon((qTm)^{1/2}N^{-1+\varepsilon})$, as $N$ goes to infinity, where $\delta_\nu(p)=0$ if $\nu=1$,
$\delta_\nu(p)=\frac1{p-p^{-1}}$ if $\nu=2$ and $\delta_\nu(p)=\frac1p$ otherwise''. For fixed $q$ and $\alpha$ the main term
dominates only for $m\le N^{1-\varepsilon}$.
\item Ellenberg \cite{Ell}, after his Corollary~2, asks ``what is the shape of the region in $m, q, N$-space for which $(a_m,L_\chi)$ is
close to $4\pi\chi(m)$?'', and remarks that a further expansion of his argument gives closeness for all $m$ up to some
$f_\delta(N)\gg N^{1-\epsilon}$.
\item \emph{Composite and prime-power indices.} Sawin and Sutherland write \cite[\S1.4]{SS}: ``Replacing the average over $a_p$ with an
average over $a_n$ for all $n$ in the study of murmurations was suggested by Bober, Booker, Lee, and Lowry-Duda [7, Remark (9) on p.~5]'',
where their [7] is \cite{BBLLD}; their Proposition~7 gives the limit of $\mathbb E_{\{E:H(E)\le X\}}\bigl[N(E)^{-1}\sum_n
W(n/N(E))\epsilon(E)a_n(E)\bigr]$ for smooth compactly supported $W$, for curves ordered by naive height $H(E)$. For squares of primes,
Kundu and M\"uller find for weight-$2$ newforms of squarefree level that ``the shape of the murmuration density is the same''
\cite[Abstract]{KM}, and add: ``It will also be straightforward to extend this idea to the case where $P$ is replaced by products of two
(or more) primes'' \cite[\S1]{KM}. Hua and Mok prove, for weight-$2$ newforms of squarefree level and every fixed $k$, that the
root-number-weighted averages at prime powers $P^k$, placed at $P^k/X$, follow the profile of the primes, under short-window and growth
conditions \cite[Theorem~1.2]{HM}. Section~\ref{sec:reduction}(b) concerns products $qp$ of a fixed small prime $q$ and a prime $p\asymp
N/q$ in the family of rank-$0$ elliptic curves ordered by conductor, which these results do not cover.
\item \emph{Approximate functional equations.} Cowan ``demonstrate[s] how murmurations arise from the averaging of approximate functional
equations'' \cite[Abstract]{Cow}; Remark~\ref{rem:series} is an exact identity of this kind for the covariance with $L(E,1)$.
\end{itemize}
For central-value-weighted averages over all newforms, asymptotics are available for $p\le N^{1-\varepsilon}$. For the covariance with
$L(E,1)$ in our family, Remark~\ref{rem:series} rewrites the dependence at $p\asymp N$ in terms of averages of $a_{pn}$ reweighted by
$w_n(E)$, a weight that depends on $N_E$, and the profile at $y$ itself \eqref{eq:split}; we know of no derivation, for elliptic curves ordered by conductor, of these averages,
of the regularity of Section~\ref{sec:reduction}(b) or of the quartile comparison of Section~\ref{sec:L}. We read, in full or in the
sections relevant here: arXiv:2204.10140, 2306.10425, 2307.00256, 2310.07681, 2408.12723, 2412.18576, 2504.09944, 2504.12295, 2505.01093,
2506.01640, 2603.09680, 2608.01337, 1605.02440 and math/0311131; Sections~1 and~3 of 1806.05869, Sections~1 and~5 of math/0510113 and
Section~1 of 2507.00738 and 2609.00777; the abstracts and
introductions of 0709.4668, 0810.4726 and 2109.07244; the abstract of 2310.07746; Sections~2.10--2.11 of \cite{Cre}; and the abstracts of
the $99$ arXiv entries on murmurations up to 2026-09-29. We did not search MathSciNet, zbMATH or Google Scholar directly.

\emph{Companion work.} The function-field family
$y^2=x^3+x+D(t)$ over $\mathbb F_q(t)$ is treated in arXiv:2603.13802v2 \cite{W1}. The main results of the companion paper
arXiv:2603.22807 are withdrawn in its second version \cite{W3}.

\subsection{Open questions}
\begin{enumerate}[leftmargin=1.8em]
\item Can the relation of Section~\ref{sec:reduction}(b) between $\Cov(\ap,a_q\mid y)$ and the rank-$0$ profile at $qy$ be
derived for curves ordered by conductor?
\item Does the first sign change of $\Cov(\ap,L(E,1)\mid y)$ sit at a fixed $y$ as $N\to\infty$? For $N<10^5$ it decreases
from $y\approx0.0068$ to $0.0045$ across conductor windows (Section~\ref{sec:L}). The positive region below it lies at $p$ of order
$\sqrt N$ ($p^*/\sqrt{N_{\mathrm{med}}}$ from $0.76$ to $1.41$ for the family part), and across the windows neither $p/N$ nor $p/\sqrt N$
is fixed. Five windows whose pairs-weighted median conductors differ by a factor $6.6$ cannot decide between these scalings, or whether either holds in the limit.
\item Does the offset $\beta$ of Section~\ref{sec:ss} tend to $0$ as $N$ grows?
\item Can the amplitude $c_q$ of Section~\ref{sec:reduction}(b) be predicted from $\ap a_q=a_{pq}$ ($p\ne q$) and a model for
averages of $a_{pq}$ at $p\asymp N/q$?
\item Do the all-$n$ densities of Sawin and Sutherland \cite[Proposition~7]{SS} or of Bober, Booker, Lee and Lowry-Duda
\cite{BBLLD}, evaluated at $n=qp$, reproduce Figure~\ref{fig:dilation}?
\item The nearest computable benchmark for Section~\ref{sec:L} is the average of $a_p(f)$ weighted by
$L(\frac12,f)L(\frac12,f\otimes\chi_D)/(f,f)$ over weight-$2$ newforms of prime level, which Michel and Ramakrishnan give exactly
\cite[Theorem~1]{MR} and Martin plots at this scale \cite[\S7]{Mar}. It differs from our statistic in weight, averaging and family. We
have not computed it at matching levels; that comparison is open.
\end{enumerate}

\subsection{Corrections to the earlier versions}\label{sec:corr}
The earlier versions (v1 and v2: arXiv:2603.04604v1 and v2, March 2026) are superseded. The main items:
\begin{itemize}[leftmargin=1.5em]
\item \emph{Definitional withdrawal.} The claim that $|\Sha|$ carries information about Frobenius traces beyond the other BSD
invariants taken together is withdrawn; read one invariant at a time, it is not a statement about $|\Sha|$ (Section~\ref{sec:def}).
\item \emph{Data provenance.} The BSD invariants are Cremona's stored values (the SageMath interface named in v1,
\texttt{LargeCremonaDatabase}, reads these tables). v1 listed the root number among the stored quantities; it is set to
$(-1)^{\mathrm{rank}}$, and PARI's \texttt{ellrootno} agrees on all $437{,}226$ class representatives (Section~\ref{sec:data}).
\item \emph{Isogeny remark.} The remark added in v2 that $|\Sha|$ is constant on isogeny classes is false (Remark~\ref{rem:isogeny}).
\item \emph{Scaling.} The claim that the amplitude of the $|\Sha|$ modulation decays as $N^{-1/4}$, ``a hallmark of genuine arithmetic
signals'', is withdrawn: a permutation null computed on the same conductor windows decays as well.
\item \emph{Zeros.} The claim that the $\Shaan\ge4$ curves have displaced low-lying zeros that connect to the difference in profiles is
withdrawn, by Section~\ref{sec:def} alone: by \eqref{eq:loglin} the indicator of $\Shaan\ge4$ is a function of $L(E,1)$, $\Omeg$, $\Tam$
and $\Tors$. The zero data reproduce; a joint test of the indicator given these and further covariates ($p=0.32$) checks only the functional
form.
\end{itemize}
Every other statement of v1 that was replaced is listed, with its status, in \texttt{anc/CORRECTIONS\_V1.md}.

\section*{Acknowledgments}
The elliptic curve data are from John Cremona's tables; the traces were computed with PARI/GP. Claude (Anthropic) was used to assist with the computations, the mathematical arguments and the writing.

\appendix
\sloppy
\section{Computation and reproducibility}\label{app:comp}

\emph{Pipelines.} (i) Cremona's \texttt{allbsd} files parsed to one table ($3{,}064{,}705$ rows; the stored values are used verbatim).
(ii) The stored trace matrix ($657{,}396\times500$, $p\le3571$; ecdata row order), checked against PARI \texttt{ellap} on $5{,}000$ random
curves with $0$ mismatches in the audit that preceded this revision. (iii) \texttt{data\_pN}: \texttt{build\_reps.py} selects curve~$1$ of each class;
$16$ PARI/GP chunks compute $\ap$ at the $9{,}592$ primes $p\le10^5$ (\texttt{run\_all.sh}, \texttt{parse\_chunk.py}); \texttt{validate.py}
compares the columns $p\le3571$ with the stored matrix; \texttt{pn\_01\_binsums.py} accumulates, per class and bin, the sums of $\ap$, $\ap^2$,
$\ap/\sqrt p$, $\ap^2/p$ and the numbers of good and bad primes ($2{,}243{,}660{,}100$ good and $1{,}581{,}248$ bad (class, prime) pairs
with $p\le N$); \texttt{pn\_00\_bruteforce\_check.py} checks the bin sums; \texttt{pn\_02\_analysis.py} produces the result tables
(bootstrap over classes, seed $20261001$, $B=1000$; $B=300$ for the composition tests, except the type-fixed column of
Table~\ref{tab:weights}, recomputed with $B=1000$ by \texttt{s48\_typefixed\_boot.py}). (iv) The scripts of Tables~\ref{tab:r2}--\ref{tab:weights}, of Section~\ref{sec:robust}, of the Sawin--Sutherland
comparison (\texttt{ss\_01}--\texttt{ss\_05}) and of the figures, listed in \texttt{anc/README.md}. (v) The scripts of this version, each a
single process reading the traces and bin sums only for the rows and columns it needs (peak memory at most $3.7$~GB):
\texttt{s01}, \texttt{s23}, \texttt{s41}, \texttt{s46} (Sections~\ref{sec:intro}--\ref{sec:data}); \texttt{s20}, \texttt{s33},
\texttt{s47}, \texttt{s60} (Section~\ref{sec:L}; \texttt{s33} and \texttt{s60} from the trace chunks); \texttt{s20}, \texttt{s24},
\texttt{s30}, \texttt{s36}, \texttt{s42} and \texttt{s70}--\texttt{s72} (Section~\ref{sec:reduction}; \texttt{s70} first reproduces the
point-evaluation outputs of \texttt{s20}; \texttt{s32} and the plateaus of \texttt{s42} are not used in the text); \texttt{s45}, \texttt{s61},
\texttt{pn\_check\_condbins.py} (Section~\ref{sec:robust}); \texttt{s31}, \texttt{s43}, \texttt{s48} (Section~\ref{sec:weights});
\texttt{s22}, \texttt{s34}, \texttt{s44}, \texttt{s62} (Section~\ref{sec:ss}); \texttt{s00} reads other numbers from their source cells.
\texttt{anc/README.md} describes each script and maps every number, table and figure to its script and input.

\emph{Validation.} New traces against the stored matrix: $218{,}613{,}000$ entries, $0$ mismatches. Brute-force checks on the full data:
the $K=40$ root-number-weighted sum and the rank-$0$ and rank-$1$ profiles recomputed from the raw traces with a different exact binning
agree with the bin sums to $0$ over all $2{,}243{,}660{,}100$ good pairs; for the $11{,}693$ classes with $N\le3571$ (all of whose primes
$p\le N$ are among the $500$ stored ones) the bin sums recomputed from the stored matrix agree exactly in the integer sums and to
$2.3\cdot10^{-13}$ in the normalised sums, on both grids. Further fixture checks are in \texttt{anc/data\_pN/README.md}. The identity \eqref{eq:loglin}
holds to $1.0\cdot10^{-14}$. For the $K=40$ profiles the ratio of bootstrap to cluster SE has median $1.000$ (range $0.95$--$1.07$).

\emph{Ancillary files} (\texttt{anc/}; \texttt{anc/README.md} lists every file and maps every table and figure to its script and
input): \texttt{CORRECTIONS\_V1.md}; \texttt{data\_pN/} (code, result CSVs, logs); \texttt{scripts/} (the scripts of (iv) and (v));
\texttt{csv/} (their outputs); \texttt{p3571\_weights/} (the weighting at $p\le3571$ on all curves).

\emph{Not checked.} The ratio $p^*/\sqrt N$ uses one representative conductor $N_{\mathrm{med}}$ per window. The traces at $p>3571$ rest on PARI \texttt{ellap} alone, apart from the classes with $N\le3571$; no second
implementation was run at large $p$. Only Cremona's BSD values are used; a PARI/GP recomputation on $1{,}600$ random curves in the
audit agreed with them to $5\cdot10^{-15}$ in $\Omeg$ and $L(E,1)$, and exactly in torsion, Tamagawa product and rank (rank: except $15$ curves of rank $\ge1$ where PARI's \texttt{ellrank} timed out). The bootstraps resample isogeny classes; clustering on the
$j$-invariant families of Section~\ref{sec:data} was checked only for Table~\ref{tab:gap40} and Figure~\ref{fig:cov}. Beyond the Bonferroni
counts stated, there is no multiple-testing adjustment across bins, weights and statistics. The standardised gap of
Section~\ref{sec:robust}, the class-level weights of Section~\ref{sec:weights}, the root mean squares of \eqref{eq:split} and the
correlations, maximisers and plateaus of Section~\ref{sec:reduction}(b) are given without standard errors (checked only by the split
sample, with point evaluation, and the two conductor windows), and the regularity of Section~\ref{sec:reduction}(b) is not derived. Apart from the
offset fits, the $\chi^2$ values of Section~\ref{sec:ss} use diagonal errors; the Hotelling calibration of the full-covariance fits assumes
Gaussian class vectors, and the comparison does not use Sawin--Sutherland's own data. The family part of Section~\ref{sec:L} truncates
the Hecke part at $pm\le2000$ and at $p\le2999$ (the omitted terms are below $10^{-16}$).
The conditioning on $L(E,1)$ was not repeated at rank $\ge1$. The search of Section~\ref{sec:lit} is limited to the sources listed there.


\begin{thebibliography}{99}

\bibitem{BBFHS} A.~Babei, B.~S.~Banwait, AJ~Fong, X.~Huang and D.~Singh, \emph{Machine learning approaches to the Shafarevich--Tate
group of elliptic curves}, arXiv:2412.18576v2, 2024.

\bibitem{BBLLD} J.~Bober, A.~R.~Booker, M.~Lee and D.~Lowry-Duda, \emph{Murmurations of modular forms in the weight aspect},
Algebra Number Theory \textbf{20} (2026), no.~2, 299--332; arXiv:2310.07746.

\bibitem{Bet} S.~Bettin, \emph{The first moment of twisted Hecke $L$-functions with unbounded shifts}, Funct.\ Approx.\ Comment.\ Math.\
\textbf{57} (2017), no.~1, 105--114; arXiv:1605.02440.

\bibitem{Cow} A.~Cowan, \emph{On the mean value of $\mathrm{GL}_1$ and $\mathrm{GL}_2$ $L$-functions, with applications to
murmurations}, arXiv:2504.09944v3, 2025.

\bibitem{Cre} J.~E.~Cremona, \emph{Algorithms for Modular Elliptic Curves}, 2nd ed., Cambridge University Press, 1997; online edition
\url{https://johncremona.github.io/book/fulltext/index.html}.

\bibitem{ecdata} J.~E.~Cremona, \emph{Elliptic curve data} (\texttt{ecdata}), \url{https://github.com/JohnCremona/ecdata}; the
\texttt{allbsd} and \texttt{allcurves} files for $N\le 499{,}998$.

\bibitem{Ell} J.~S.~Ellenberg, \emph{On the error term in Duke's estimate for the average special value of $L$-functions}, Canad.\
Math.\ Bull.\ \textbf{48} (2005), no.~4, 535--546; arXiv:math/0311131v2.

\bibitem{Faz} A.~Fazzari, \emph{A weighted one-level density of families of $L$-functions}, Algebra Number Theory \textbf{18} (2024),
no.~1, 87--132; arXiv:2109.07244.

\bibitem{FW} B.~Feigon and D.~Whitehouse, \emph{Averages of central $L$-values of Hilbert modular forms with an application to
subconvexity}, Duke Math.\ J.\ \textbf{149} (2009), no.~2, 347--410; arXiv:0810.4726.

\bibitem{HLOP} Y.-H.~He, K.-H.~Lee, T.~Oliver and A.~Pozdnyakov, \emph{Murmurations of elliptic curves}, Exp.\ Math.\ \textbf{34}
(2025), no.~3, 528--540; arXiv:2204.10140v2.

\bibitem{HM} S.~Hua and C.~P.~Mok, \emph{Universal murmuration and Hecke augmentation}, arXiv:2609.00777v1, 2026.

\bibitem{KR} A.~Knightly and C.~Reno, \emph{Weighted distribution of low-lying zeros of $\mathrm{GL}(2)$ $L$-functions}, Canad.\ J.\
Math.\ \textbf{71} (2019), no.~1, 153--182; arXiv:1806.05869v2.

\bibitem{KM} D.~Kundu and K.~M\"uller, \emph{Murmurations of modular forms and $p$-power coefficients}, Math.\ Proc.\ Cambridge
Philos.\ Soc.\ \textbf{180} (2026), no.~1, 65--92; arXiv:2507.00738v1.

\bibitem{Mar} K.~Martin, \emph{Variations on murmurations}, arXiv:2505.01093v1, 2025.

\bibitem{MR} P.~Michel and D.~Ramakrishnan, \emph{Consequences of the Gross--Zagier formulae: stability of average $L$-values,
subconvexity, and non-vanishing mod $p$}, in \emph{Number Theory, Analysis and Geometry}, Springer, 2012, 437--459; arXiv:0709.4668.

\bibitem{PARI} The PARI~Group, \emph{PARI/GP version 2.17}, Univ.\ Bordeaux, \url{https://pari.math.u-bordeaux.fr/}.

\bibitem{RR} D.~Ramakrishnan and J.~Rogawski, \emph{Average values of modular $L$-series via the relative trace formula}, Pure Appl.\
Math.\ Q.\ \textbf{1} (2005), no.~4, 701--735; arXiv:math/0510113.

\bibitem{RRe} D.~Ramakrishnan and J.~Rogawski, \emph{Erratum: Average values of modular $L$-series via the relative trace formula},
Pure Appl.\ Math.\ Q.\ \textbf{5} (2009), no.~4, 1469; doi:10.4310/PAMQ.2009.v5.n4.a10.

\bibitem{SS} W.~Sawin and A.~V.~Sutherland, \emph{Murmurations for elliptic curves ordered by height}, arXiv:2504.12295v2, 2025.

\bibitem{W1} D.~Wachs, \emph{Frobenius in $W(E_8)$ for $y^2=x^3+x+D(t)$ over $\mathbb F_q(t)$: ranks, Tate--Shafarevich groups and
trace statistics}, arXiv:2603.13802v2, 2026.

\bibitem{W2v1} D.~Wachs, \emph{BSD invariants and murmurations of elliptic curves}, arXiv:2603.04604v1 and v2, 2026 (the earlier versions
of this paper).

\bibitem{W3} D.~Wachs, \emph{Murmurations, periods, and local factors}, arXiv:2603.22807v2, 2026.

\bibitem{Zub} N.~Zubrilina, \emph{Murmurations}, Invent.\ Math.\ \textbf{241} (2025), no.~3, 627--680; arXiv:2310.07681.

\end{thebibliography}
\end{document}